\documentclass[a4paper]{article}

\usepackage{geometry}
\usepackage[utf8]{inputenc}
\usepackage[english
]{babel}
\usepackage{microtype}
\usepackage{amsmath, amsfonts}
\usepackage{amsthm}
\usepackage[table]{xcolor}  
\usepackage{tikz}
\usetikzlibrary{shapes}
\usetikzlibrary{shapes.misc}
\usetikzlibrary{snakes}
\usetikzlibrary{trees}
\usepackage{graphicx}
\graphicspath{{fig/}}
\usepackage{import}
\usepackage{ragged2e}
\usepackage{array}
\usepackage{rotating}
\usepackage{paralist}
\usepackage[square,comma]{natbib}
\usepackage{float}
\usepackage{listings}
\usepackage{multirow}
\usepackage{multicol}
\usepackage{booktabs}
\usepackage[perpage]{footmisc}

\usepackage[pdfpagelabels,plainpages=false,pdfborder={0 0 0}]{hyperref}
\usepackage{algorithm}

\usepackage{subfigure}

\newtheorem{lem}{Lemma}
\newtheorem{theorem}[lem]{Theorem}

\newtheorem{coro}[lem]{Corollary}

\newcommand{\bydef}{\stackrel{\rm{def}}{=}}

\newcommand{\toN}{^{N}}
\newcommand{\E}{\mathbb{E}}
\newcommand{\N}{\mathbb{N}}
\newcommand{\R}{\mathbb{R}}
\newcommand{\A}{\mathcal{A}}
\renewcommand{\Pr}{\mathbb{P}}
\newcommand{\PP}{\mathcal{P}}
\renewcommand{\SS}{\mathcal{S}}
\newcommand{\re}{\mathrm{r}}

\renewcommand{\time}{\times}
\newcommand{\CC}{\mathrm{C}}
\newcommand{\ccc}{c}
\newcommand{\MM}{\mathrm{M}}

\newcommand{\MRn}{\Pr_N(\SS){\time}\R^d}

\newcommand\tolaw{\xrightarrow{\mathcal{L}}}
\newcommand\toas{\xrightarrow{a.s}}

\newcommand\norminf[1]{\|#1\|_\infty}
\newcommand\normun[1]{\|#1\|_1}

\title{A Mean Field Approach for Optimization in Particles Systems and
  Applications}
\author{Nicolas Gast \and Bruno Gaujal}

\usepackage{RR}

\RRdate{March 2009}

\RRversion{2}
\RRdater{June 2009}

\RRauthor{Nicolas Gast \and Bruno Gaujal}
\authorhead{N.Gast \& B. Gaujal}

\RRtitle{Une approche champ moyen pour l'optimisation dans les systèmes de
  particules et ses applications}
\RRetitle{A Mean Field Approach for Optimization in Particles Systems and
  Applications}
\titlehead{A Mean Field Approach for Optimization in Particles Systems and
  Applications}
\RRresume{Cet article examine le comportement limite de processus de
  décision Markovien constitués de particules indépendantes évoluant dans un
  environnement commun, lorsque le nombre de particules tend vers l'infini.
  
  Dans le cas où on s'intéresse à un coût à horizon fini ou dans le cas
  d'un coût à horizon infini avec décote, nous montrons que lorsque le
  nombre de particules devient grand, le coût optimal du système converge
  presque sûrement vers le coût optimal du système déterministe. La
  convergence vaut également pour les politiques optimales.
  
  De plus, nous donnons un aperçu de la vitesse de convergence en prouvant
  plusieurs théorèmes de la limite centrale pour le coût ainsi que l'état
  moyen du processus en donnant des formules explicites pour la variance
  des lois gaussiennes limites.

  Enfin, ce modèle est appliqué à un problème de gestionnaire de ressources
  dans des grilles de calcul. Nous donnons un algorithme explicite pour
  calculer la politique optimale de la limite puis plusieurs simulations
  avec un nombre variable de processeurs sont étudiées. Nous comparons les
  performances de la politique optimale de la limite appliquée au système
  initiale avec plusieurs politiques classiques, (telles que joindre la
  file la plus courte). Nous mesurons le gain asymptotique, ainsi que le
  seuil à partir duquel elle surpasse les politiques classiques.}

\RRabstract{  This paper investigates the limit behavior of Markov decision processes
  (MDPs) made of independent particles evolving in a common environment,
  when the number of particles goes to infinity.

  In the finite horizon case or with a discounted cost and an infinite
  horizon, we show that when the number of particles becomes large, the
  optimal cost of the system converges almost surely to the optimal cost of
  a deterministic system (the ``optimal mean field'').
  Convergence also holds for optimal policies.
  
  We further provide insights on the speed of convergence by proving
  several central limits theorems for the cost and the state of the Markov
  decision process with explicit formulas for the variance of the limit
  Gaussian laws.

  Then, our framework is applied to a brokering problem in grid computing.
  The optimal policy for the limit deterministic system is computed
  explicitly. Several simulations with growing numbers of processors are
  reported. They compare the performance of the optimal policy of the limit
  system used in the finite case with classical policies (such as Join the
  Shortest Queue) by measuring its asymptotic gain.
}

\RRmotcle{Processus de décision Markovien, Champ moyen,
  Optimisation, Systèmes de particules, Gestionnaire de ressource}
\RRkeyword{Markov Decision Processes, Mean Field, Optimization,  Particles
  System, Grid Broker}

\RRprojet{MESCAL}
\RRtheme{\THNum} 
\URRhoneAlpes 

\begin{document}
\RRNo{6877} 

\makeRR

\section{Introduction}

The general context of this paper is the optimization of the behavior
of controlled Markovian systems, namely Markov Decision Processes 
composed by a large number of particles evolving in a common
environment.

Consider a discrete time system made of $N$ particles, $N$ being large,
that evolve randomly and independently (according to a transition
probability kernel $K$).  At each step, the state of each particle changes
according to a probability kernel, depending on the environment. The
evolution of the environment only depends on the number of particles in
each state.  Furthermore, at each step, a central controller makes a
decision that changes the transition probability kernel.  The problem
addressed in this paper is to study the limit behavior of such systems when
$N$ becomes large and the speed of convergence to the limit.

Several  papers (\cite{benaim:cmf}, \cite{bordenave2007psi}) 
study the limit behavior of Markovian systems in the case of vanishing intensity (the expected number of transitions
per time slot is  $o(N)$).  In these cases, the system converges
to a  differential system in continuous time. In the case considered here,
time remains 
discrete at the limit. This requires a rather different approach to
construct the limit.

In \cite{boudec2007gmf}, discrete time systems are considered and the authors show that under certain
conditions, as $N$ grows large, a Markovian system made of $N$
particles  converges to a deterministic system. Since a Markov
decision process  can be seen as a family of Markovian  kernels, 
the class of systems studied in  \cite{boudec2007gmf} corresponds to
the case where this  family is reduced to
a unique kernel  and  no decision can be
made.
Here,  we show that under similar conditions as in  \cite{boudec2007gmf},
a Markov decision process  also converges to a deterministic one.  
More precisely, we show that the optimal costs (as well as the
corresponding states)
converge almost surely to the optimal costs (resp. the
corresponding states) of a deterministic system (the ``optimal mean field'').

On a practical point of view, this allows one  to
compute  the optimal policy in a deterministic system which can
often be done very efficiently, and then to use this policy in the
original  random system as a good approximation of the optimal policy,
which cannot be  computed efficiently because of the curse of
dimensionality.
This is illustrated by an application of our framework to optimal
brokering in computational grids.
We consider a set of multi-processor clusters (forming a computational grid,
like EGEE \cite{EGEE}) and a set of users submitting tasks to be executed. 
A central broker assigns the tasks to the clusters (where tasks are
buffered and served in a fifo order) and tries to
minimize the average processing time of all tasks.
Computing the optimal policy (solving the associated MDP) is
known to be hard \cite{Tsitsiklis}. Numerical computations can only be carried
up to a total of 10 processors and two users.
However, our approach shows that 
when the number of processors per cluster and the number of users
submitting tasks grow, the system converges to a mean field
deterministic system.
For this deterministic mean field system,  the optimal brokering policy can
be explicitly computed. 
Simulations reported in Section \ref{sec:example} show that, using this
policy over a grid with a growing number of processors, makes
performance converge to the optimal sojourn time in a deterministic
system, as expected. Also, simulations show that this deterministic
static policy
outperforms classical dynamic policies such as Join the Shortest
Queue,  as soon as
the total number of processors and users is  over 50.

In general, how good  the deterministic  approximation is and how fast convergence takes place
can also be estimated. For that, we 
provide bounds on the speed of convergence by proving  of central limit
theorem for the state of the system under the optimal policy as well
as for the cost function.

\section{Notations and definitions}

The system is composed of $N$ particles. There are $S$ possible
states for each particle, the state space is denoted by
$\mathcal{S}{=}\{1,\dots,S\}$. The state of the $n$th particle at time $t$
is denoted $X\toN_n(t)$. We assume that the particles are distinguishable
only through their state and that the dynamics of the system is homogeneous
in $N$. In other words, this means that the behavior of the system only
depends on in the proportion of particles in every state $i$.  For all
$i\in\SS$, $\big(\MM\toN_t\big)_i\bydef \sum_{n=1}\toN
\mathbf{1}_{X\toN_n(t)=i}$ is the proportion of particles in state $i$ and
we denote by $\MM\toN_t$ the vector
$(\big(\MM\toN_t\big)_1\dots\big(\MM\toN_t\big)_S)$. The set of possible
values for $\MM\toN$ is the set of probability measures $p$ on $\{1\dots
S\}$, such that $Np(i)\in\N$ for all $i\in\SS$, denoted by
$\mathcal{P}_N(S)$. For each $N$, $\mathcal{P}_N(S)$ is a finite set. When
$N$ goes to infinity, it converges to $\PP(\SS)$ the set of probability
measures on $\SS$.

The system of particles evolves depending on their common environment. We
call $\CC\in\R^d$ the context of the environment. 
Its evolution depends on the mean states of the particles $\MM\toN$, itself
at the previous time slot and the action $a_t$ chosen by the controller
(see below):
\begin{equation*}
  \CC\toN_{t+1} = g(\CC\toN_t, \MM\toN_{t+1},a_t),
\end{equation*}

where $g: \PP_N(\SS){\time}\R^d$ $ {\time} \mathcal{A} \to$ $\R^d$ is a continuous function.

\subsection{Actions and policies}

At each time $t$, the system's state is $\MM\in\PP_N(\SS)$. The decision
maker may choose an action $a$ from the set of possible actions $\A$. $\A$
is assumed to be a compact set   (finite or infinite). The action
determines how the system will evolve. For an action $a\in\A$ and an
environment $\CC\in\R^d$, we have a transition probability kernel $K(a,\CC)$ such
that the probability that a particle  goes  from state $i$ to state the
$j$ is $K_{i,j}(a,\CC)$:
\begin{equation*}
  \Pr(X\toN_n(t+1) = j | X\toN_n(t) = i, a_t = a, \CC\toN_t = \CC) = K_{i,j}(a,\CC). 
\end{equation*}

The evolutions of particles are supposed to be independent once $\CC$ is
given. Moreover, we assume that $K_{i,j}(a,\CC)$ is continuous in $a$ and
$\CC$. The assumption of independence of the users is a rather common
assumption in mean field models \cite{boudec2007gmf}. However other papers
\cite{benaim:cmf,bordenave2007psi} have shown that similar results can be
obtained using asymptotic independence only (see \cite{Graham} for results
of this type).

Here, the focus is on Markov Decision Processes theory and on the
computation of  optimal
policies. A policy $\Pi=(\Pi_1\dots\Pi_t\dots)$ specifies the decision
rules to be used at each time slot. A decision rule $\Pi_t$ is a procedure
that provides an action at time $t$.  In general, $\Pi_t$ is a random
measurable function that depends on the events
$((\MM_1,\CC_1)\dots(\MM_t,\CC_t))$ but it can be shown that when the state
space is finite and the action space is compact, then deterministic
Markovian policies (\emph{i.e.} that only depends deterministically on the current state) are
dominant,  therefore we will only focus on them
\cite{puterman1994mdp}.

\subsection{Reward functions}
To each possible state $(\MM,\CC)$ of the system at time $t$, we associate a
reward $\re_t(\MM,\CC)$. The reward is assumed to be continuous  in $\MM$ and $\CC$. This 
function can be either seen as a reward -- in that case the controller 
wants to maximize the reward --, or as a cost -- in that case the goal of
the controller is to minimize this cost. In this paper, we will focus on
two problems: finite-horizon reward and discounted reward.

In the finite-horizon case, we want to maximize the sum of the rewards over all 
time $t<T$ plus a  final reward that depends on the final state, 
$\re_T(\MM\toN_T,\CC\toN_T)$. The expected reward of the policies
$\Pi_0,\dots,\Pi_{T-1}$ is:
\begin{equation*}
  V\toN_{\Pi_0\dots \Pi_T}(\MM\toN_0,\CC\toN_0) \bydef
  \E\bigg[\sum_{t=1}^{T-1} \re_t(\MM\toN_t,\CC\toN_t)  + \re_T(\MM\toN_T,\CC\toN_T)
  \bigg], 
\end{equation*}

where the expectation is taken over all possible $(\MM\toN_t,\CC\toN_t)$ when
the actions are $\Pi_t(\MM\toN_t,\CC\toN_t)$, for all $t$.

Let $0\leq\delta<1$, the discounted reward associated to  $\delta$ and the
policy $\Pi_0\dots\Pi_t\dots$ is the quantity:
\begin{equation*}
  V\toN_{(\delta),\Pi_0\dots}(\MM\toN_0,\CC\toN_0) \bydef 
  \E\bigg[\sum_{t=1}^\infty \delta^t \re_t(\MM\toN_t,\CC\toN_t)\bigg].
\end{equation*}

Again, the expectation is taken  over  all possible $(\MM\toN_t,\CC\toN_t)$ when
the actions at time $t$ is $\Pi_t(\MM^N_t,\CC^N_t)$, for all $t$.

In both cases, the goal of the controller is to find a policy that
maximizes the expected reward:
\begin{equation*}
  {V^*}\toN(\MM\toN_0,\CC\toN_0) \bydef  \sup_{\Pi_1\dots \Pi_T}
  V\toN_{\Pi_1\dots\Pi_T}(\MM\toN_0,\CC\toN_0),
\end{equation*}
\begin{equation*}
V^{*N}_{(\delta)}(\MM\toN_0,\CC\toN_0) \bydef  \sup_{\Pi_1\dots}
  V\toN_{(\delta),\Pi_1\dots}(\MM\toN_0,\CC\toN_0).
\end{equation*}

\subsection{Summary of the assumptions}
\label{sec:assum}

Here is the list of  the assumptions under which all our results
will hold, together with some comments on their tightness and  their degree of
generality and applicability.

\newcounter{Lcount}
\begin{list}{(A\arabic{Lcount})}{\usecounter{Lcount}
    \setlength{\itemindent}{0.4cm}
    \setlength{\leftmargin}{0.5\leftmargin}
  }
\item \label{assum:indep} \textbf{Independence of the users, Markov system} --
  If at time $t$ if the environment is $\CC$ and the action is $a$, then
  the behavior of each particle is independent of other particles and
  its evolution is Markovian with a kernel $K(a,\CC)$. 

\item \label{assum:action} \textbf{Compact action set} -- The set of action $\A$ is
  compact. 
\item \label{assum:conti} \textbf{Continuity of $K,g,\re$} -- the mappings
  $(\CC,a)\mapsto K(a,\CC)$, $(\CC,\MM,a)\mapsto g(\CC,\MM,a)$ and $(\MM,\CC)\mapsto
  \re_t(\MM,\CC)$ are continuous deterministic functions,
  uniformly continuous in $a$. 
\item \label{assum:init_as} \textbf{Almost sure initial state} -- Almost
surely, the initial measure $\MM\toN_0,\CC\toN_0$ converges to a
  deterministic value $m_0,\ccc_0$.  Moreover, there exists $B<\infty$ such
  that almost surely $\norminf{\CC\toN_0}\leq B$ where $\norminf{C}=\sup_i
  |C_i|$. 
\end{list}


 To
simplify the notations, we choose the functions $\CC$ and $g$ not to depend on
time. However as the proofs will be done for each time step, they also
hold if the functions are time-dependent (in the finite horizon
case). 

Also,  $K,g$ and $\re$ do not to depend on $N$, while this is the
  case in most practical cases. Adding a uniform continuity assumption on
  these functions for all $N$ will make all the proofs work the
  same. 



Here are some comments on the uniform bound $B$ on the initial
  condition (A4). In fact, as $\CC\toN_0$ converges almost surely,
  $\CC\toN_0$ is almost surely bounded. Here we had a bound $B$ which is
  uniform on all events in order to be sure that the variable $\CC\toN_0$
  is dominated by an integrable function. As $g$ is continuous and the sets
  $\A$ and $\PP(\SS)$ are compact, this shows that for all
  $t$, 
there exists $B_t<\infty$ such that  
\begin{equation}
  \norminf{\CC\toN_t} \leq B_t. \label{eq:unif_bounded}
\end{equation}

Finally, in many cases the rewards also depend  on the
action. This is not the case here, at a small loss of generality.

\section{Convergence results and optimal policy}

In the case where there is no control, one can adapt the results proved in
\cite{boudec2007gmf} to show that when $N$ goes to infinity, the system
converges almost surely to a deterministic one. In our case, this means
that if the actions are fixed, the system converges.

For any fixed action $a$ and any value $\MM\in\PP_N(\SS)$, we define the
random variable $\Phi\toN_a(\MM,\CC)$ that corresponds to the state of the
system $\MM',\CC'$ after one iteration started from $\MM,\CC$. For
$m\in\Pr(\SS)$, we define $\Phi_a(m,\ccc)$ the (deterministic) value
corresponding to one iteration of the mean field system:
$\Phi_a(m_t,\ccc_t)=(m_{t+1},\ccc_{t+1})$ where 
\begin{eqnarray*}
m_{t+1} &=& m_t.K(a,\ccc_t)\\
  \ccc_{t+1} &=& g(m_{t+1},\ccc_t).
\end{eqnarray*}

We call $\Phi\toN_{a_0\dots a_{T-1}}$ (resp. $\Phi_{a_0\dots a_{T-1}}$) the
compositions of $\Phi\toN_{a_0},\dots,\Phi\toN_{a_{T-1}}$ (resp. of
$\Phi_{a_0}\dots\Phi_{a_{T-1}}$).

In \cite{boudec2007gmf}, the system is homogeneous in time. However,  the
proofs are done for each step time and the results still hold without
time homogeneity. With our
notations, theorem 4.1 of \cite{boudec2007gmf} says that if the actions are
$a_0\dots a_{T-1}$, and if the initial state converges almost surely, then
the system of size $N$ converges almost surely.
\begin{theorem}[Mean Field Limit, th. 4.1 of \cite{boudec2007gmf}] 
  \label{th:convergence_lebou}
  Under assumptions (A1,A3,A4), if the controller takes the actions
  $a_t$ at time $t$, then for any fixed $T$: 
  \begin{equation*}
    (\MM\toN_t,\CC\toN_t) \toas \Phi_{a_0\dots a_{T-1}}(m_0,\ccc_0).
  \end{equation*}
\end{theorem}

In the following, we will first show that if we fix the actions, the total
reward of the system converges when $N$ grows, then we will show
that the optimal reward also converges.

\subsection{Finite horizon model}
\label{sec:finite-horizon}

In this section, the horizon  $T$ is fixed, the infinite horizon case will be
treated in Section \ref{sec:infinite-discouted}.
Using the same notation and hypothesis as in
Theorem~\ref{th:convergence_lebou}, we define the reward of the deterministic
system starting at $m_0,\ccc_0$ under the actions $a_0,\dots,a_{t-1}$: 
\begin{equation*}
  v_{a_0\dots a_{t-1}}(m_0,\ccc_0) = \sum_{t=1}^T
  \re_t(\Phi_{a_0\dots a_{t-1}}(m_0,\ccc_0)).
\end{equation*}

For any $t$, if the action taken at instant $t$ is
fixed  equal to $a_t$, then $(\MM\toN_t,\CC\toN_t)$ converges almost surely to
$(m_t,\ccc_t)$. Since the reward at time $t$ is continuous, this means that the finite-horizon
expected reward converges as $N$ grows large: 

\begin{lem}[\label{lem:converg_reward}Convergence of the reward]
  Under assumptions (A1,A3,A4), if the controller takes  actions
  $a_0$ $\dots$ $a_{T-1}$, the finite-horizon expected reward of the stochastic
  system converges to the finite-horizon reward of the deterministic
  system:
  \begin{equation*}
    \lim_{N\to\infty} V\toN_{a_0\dots a_{t-1}}(\MM\toN_0,\CC\toN_0) =
    v_{a_0,\dots,a_{t-1}}(m_0,\ccc_0) \quad \mathrm{a.s}. 
  \end{equation*}
\end{lem}

\begin{proof}
  For all $t$, $(\MM\toN_t,\CC\toN_t)$ converges almost surely to
  $(m_t,\ccc_t)$. Since the reward at time $t$ is continuous in $(\MM,\CC)$,
  then $\re_t(\MM\toN_t,\CC\toN_t)\toas \re_t(m_t,\ccc_t)$. Moreover, as
  $(\MM,\CC)$ are bounded (see Equation \eqref{eq:unif_bounded}), 
  the dominated convergence theorem shows that 
  $\E[\re_t(\MM\toN_t,\CC\toN_t)]$ goes to $\re_t(m_t,\ccc_t)$ which
  concludes the demonstration. 
\end{proof}

Now, let us consider the problem of  convergence of the reward
under the optimal strategy  of the controller.
First, it should be clear that the optimal strategy exists for the
limit system.  Indeed, the limit system being  deterministic, starting
at state $(m_0,\ccc_0)$, one only needs to know the actions to take for all
$(m_t,\ccc_t)$ to compute the reward. The  optimal policy is
deterministic and 
$v^*_T(m_0,\ccc_0) \bydef  \sup_{a_0 \dots a_{T-1}} %
\{v_{a_0 \dots
  a_{T-1}}(m_0,\ccc_0)\}$. Since the action set is compact, this supremum is
a maximum: there  exist $a^*_0\dots a^*_{T-1}$ such that $v^*_T(m_0,\ccc_0)
= v_{a^*_0 \dots a^*_{T-1}}(m_0,\ccc_0)$. In fact,  in many cases
there are  more than one optimal action sequence. In the following, $a^*_0\dots
a^*_{T-1}$ is one of them, and will be called  the sequence of  {\it optimal limit actions}.

\begin{theorem}[\label{th:converg_opti}Convergence of the optimal reward]
  Under assumptions (A1,A2,A3,A4), as $N$ goes to infinity, the optimal
  reward of the stochastic system converges to the optimal reward of the
  deterministic limit system: almost surely,
  \begin{equation*}
    \lim_{N\to\infty}
    V^{*N}_T  (M_0^N,C_0^N)~{=} \lim_{N\to\infty}V\toN_{a^*_0\dots
      a^*_{T-1}}  (M_0^N,C_0^N)=  v^*_T(m_0,c_0)
  \end{equation*}
\end{theorem}

In words, this theorem says that, at the limit, the reward of the
optimal policy under full information $ {V^*}\toN_T
(M_0^N,C_0^N)$   is the same as 
the  reward obtained
when the optimal  limit actions   $(a^*_0\dots a^*_{T-1})$ are used
in the original system, 
both being equal to  the optimal reward of the limit deterministic
  system,  $v^*_T(m_0,c_0)$.

\begin{proof}
  For all $N$ and $0\leq t\leq T$ and $(\MM,\CC)\in \MRn$, let us define by
  induction on $t$ the function $V^{*N}_{t\dots T}$:
  \begin{equation}
      \begin{array}{l}
        V^{*N}_{T\dots T}(\MM,\CC) = \re_T(\MM,\CC)\\
        V^{*N}_{t\dots T}(\MM,\CC) {=} \re_t(\MM,\CC){+}\displaystyle\sup_{a\in\A}
        \E_{\MM,\CC}[V^{*N}_{t+1\dots T}(\Phi\toN_a(\MM,\CC))].
      \end{array}
    \label{eq:V_*_1...t}
  \end{equation}
  where the expectation $\E_{\MM,\CC}[\cdot]$ is taken over  all possible values of
  $\Phi\toN_a(\MM,\CC)$ given $(\MM,\CC)$.  Also notice that 
  $V^{*N}_{t\dots T}(\MM,\CC)$ is the maximal expected reward
  between time $t$ and time $T$ starting in $(\MM,\CC)$ and therefore
  $V^{*N}_{0\dots T} = V^{*N}_T$.

  Let us also define for the limit system, $v^*_{t\dots T}$ similarly (by removing the
  expectation):
  \begin{equation}
    \begin{array}{lll}
    v^*_{T\dots T}(m,c) &=& \re_T(m,c)\\
    v^*_{t\dots T}(m,c) &=& \re_t(m,c)+\displaystyle\sup_{a\in\A}
    \Big[v^*_{t+1\dots T}\big(\Phi_a(m,c)\big)\Big],
  \end{array}  \label{eq:v_*_1...t}
  \end{equation}
  and let $\Pi^*_t(m,c)$ be an action that maximize the $\sup$ in the
  previous equation (it exists because of (A2): $\A$ is compact). 
  
  We will show by induction on $t<T$ that $V^{*N}_{t\dots T}(\cdot,\cdot)$
  is continuous (note that since $\MM\in\PP^N(\SS)$ is discrete the
  continuity in $\MM$ is trivial) and that we can define an optimal policy
  $\Pi^{*N}_t(\MM,\CC)$, such
  that:
  \begin{equation}
    V^{*N}_{t\dots T}(\MM,\CC) {=} \re_t(\MM,\CC){+}
    \E\big[V^{*N}_{t+1\dots T}(\Phi\toN_{{\Pi^*_t}\toN(\MM,\CC)}(\MM,\CC))\big]. 
  \end{equation}

  For $t=T$, the assumption holds by the continuity of $\re$ (A3). 
  
  Let us assume that it holds for $t+1\leq T$. By assumption (A3), the
  mapping $g$ and the kernel $K$ are continuous in $a$ thus if
  $\{a(k)\}_{k\in \N}$ is a
  sequence of action converging to $a$, $\Phi\toN_{a(k)}$ converges (in law)
  to $\Phi\toN_a$. As $V^{*N}_{t+1\cdots T}$ is continuous,
  $a\mapsto\E[V^{*N}_{t+1\dots T}(\Phi\toN_a(\MM,\CC))]$ is continuous.
  Using 
  this continuity and the compacity of $\A$, the optimal action
  $\Pi^{*N}_t(\MM,\CC)\in\A$  exists.  The functions $\re$, $g$, $K$ are
  uniformly continuous in $a$, therefore the convergence of the continuity
  of the function $a\mapsto\sup_a\E[V^{*N}_{t+1\dots
    T}(\Phi\toN_a(\MM,\CC))]$ is uniform in $M,R$. This shows that
  $(M,R)\mapsto \sup_a\E[V^{*N}_{t+1\dots T}(\Phi\toN_a(\MM,\CC))]$ is
  continuous and  the property for all $t$ is proved.

  Let us now prove by induction on $t$ that for all sequences
  $(\MM\toN,\CC\toN)$ converging almost surely to $(m,\ccc)$, 
  $v^{*N}_{t\dots T}(\MM\toN,\CC\toN) $ $\toas v^*_{t\dots T}(m,\ccc)$. This is
  clearly true for $t{=}T$. Assume that it holds for some $t{+}1{\leq}T$
  and let 
  us call $a^*_{t}\dots a^*_{T-1}$ a sequence of optimal actions for the
  deterministic limit. Lemma~\ref{lem:converg_reward} shows that
  $V\toN_{a^*_{t}\dots a^*_{T-1}}(\MM\toN,\CC\toN) \toas v_{a^*_{t}\dots
    a^*_{T-1}}(m,\ccc) = v^*_{t\dots T}(m,\ccc)$. In particular, this shows
  the second inequality (which holds a.s.) of the following equation:
  \begin{equation}
    \begin{array}{r}
      \liminf V^{*N}_{t\dots T}(\MM\toN,\CC\toN) \geq \liminf
      V\toN_{a^*_{t}\dots a^*_{T-1}}(\MM\toN,\CC\toN) \\
      = v^*_{t\dots  T}(m,\ccc).
    \end{array}
  \end{equation}

  Let $a^{*N}$ be a sequence of actions maximizing the expectation in
  \eqref{eq:V_*_1...t}. As $\A$ is compact, there exists a subsequence
  $a^{*\psi(N)}$ converging to a value $a$. Again by
  lemma~\ref{lem:converg_reward}, the $\limsup$ of 
  $\re(\MM^{\psi(N)},\CC^{\psi(N)})+\E[V^{*{\psi(N)}}_{t+1\dots
    T}(\Phi^{\psi(N)}_a(\MM^{\psi(N)},\CC^{\psi(N)}))]$ converges a.s. to 
  $\re(m,\ccc)+v^*_{t+1}(\Phi_a(m,\ccc)) \leq v^*_{t\dots T}(m,\ccc)$. Using
  both inequalities, this shows that $V^{*{\psi(N)}}_{t\dots
    T}(\MM^{\psi(N)},\CC^{\psi(N)}) \toas v^*(m,\ccc)$.
  
  To conclude the proof, remark that since the limit system
  is deterministic and takes the values $(m_0,\ccc_0),\dots,(m_t,\ccc_t)$, fixing
  the policy at time $t$ to the action $a^*_t\bydef\Pi^*(m_t,\ccc_t)$ achieves
  the optimal reward.
\end{proof}

This result has several practical consequences.
Recall that  the limit actions $a_0^*\dots a_{T-1}^*$ is  a sequence of optimal actions in the
limit  case, \emph{i.e.} such that $v_{a^*_0\dots
  a^*_{t-1}}(m,\ccc)=v^*_T(m,\ccc)$. 
This result  proves that in the limit case, the
optimal policy does  not depend on the state of the system. This also shows that
incomplete information policies are as good as complete information
policies.  
However, the state $(\MM\toN_t,\CC\toN_t)$ is not deterministic and on
one trajectory of the system,  it could be quite far from its
deterministic limit 
$(m_t,\ccc_t)$. In the proof of proposition~\ref{lem:converg_reward},
we  also defined the policy $\Pi^*_t(\MM\toN_t,\CC\toN_t)$  which is optimal for the deterministic  
system starting at time $t$ in state $m_t,r_t$. The least we can say
is that this  
strategy is also asymptotically optimal, that is: 
\begin{equation*}
  \lim_{N\to\infty}  V\toN_{\Pi^*_0\dots \Pi^*_T}(\MM,\CC) = \lim_{N\to\infty}
  V\toN_{a^*_0\dots a^*_T}(\MM,\CC).
\end{equation*}
In practical situations, using this policy will decrease the risk of being
far from the optimal state. On the other hand, using this policy has some
drawbacks. The first one is that the complexity of computing the optimal
policy for all states can be much larger than the complexity of computing
$a^*_0\dots a^*_{T-1}$.  An other one is that the system becomes very
sensitive to random perturbations: the policy $\Pi^*$ is not necessarily
continuous and may not have a  limit. 
In Section \ref{sec:example}, a comparison between the performances of
$a^*_0\dots a^*_{T-1}$ and 
$\Pi^*_0\dots \Pi^*_{T-1} $  is provided over an example.

\subsection{\label{sec:clt}Central Limit  Theorems}

In this part we prove central limit theorems for 
interacting particles. This result  provides estimates  on the speed of
convergence to the mean field limit.  This section contains two main
results: 

The first one is that when the control action sequence is fixed,   the
gap to the mean field limit  decreases as
the inverse square root of the number of particles.
The second result states that the gap between the optimal reward for the finite
system and the optimal reward for the limit system also decreases  as fast
as $1/\sqrt{N}$. These properties are
formalized in theorems  \ref{th:clt_particles} and \ref{th:clt_reward} respectively.

To prove these results, we will need additional assumptions (A4-bis) and
(A5) or (A5-bis).

\begin{list}{}{    
    \setlength{\itemindent}{0.4cm}
    \setlength{\leftmargin}{0.5\leftmargin}
}
\item[(A\ref{assum:init_as}-bis)] \textbf{Initial Gaussian variable} -- There exists a
  Gaussian vector $G_0$ of mean $0$ with covariance $\Gamma_0$ such that
  the  vector $\sqrt{N}((\MM\toN_0,\CC\toN_0){-}(m_0,\ccc_0))$ (with $S{+}d$ components)
  converges in law to $G_0$. (This is denoted as
  $\sqrt{N}((\MM\toN_0,\CC\toN_0)-(m_0,\ccc_0))\tolaw G_0$). This
  assumption also includes (A4), {\it i.e.} almost sure convergence of the
  initial state. 
\item[(A5)] \textbf{Continuous differentiability} -- For all $t$ and
  all $i,j\in\SS$, all functions $g$, $K_{ij}$ and $\re_t$ are continuously
  differentiable.
\item[(A5-bis)] \textbf{Differentiability in $a_0\dots a_{T-1}$} -- Let
  $(m_t,\ccc_t)$ be the deterministic limit of the system if the
  controller takes the actions $a_0\dots a_{T-1}$ then for all $i,j\in\SS$, the
  functions $g$, $K_{ij}$ and $\re_t$ are differentiable in the points
  $(m_t,\ccc_t)$.
\end{list}

These assumptions  are slightly stronger than (A3) and (A4) but  remain very natural.
(A4-bis) is clearly necessary for  Theorems  \ref{th:clt_particles} and
\ref{th:clt_reward} to
hold. The differentiability condition implies that if the gap between
$\MM_t$ and $m_t$ is of order $1/\sqrt{N}$, it remains of the same order at
time $t+1$. For Theorem \ref{th:clt_particles}, (A5-bis) is
necessary but can be replaced by a Lipschitz continuity condition for
Theorem \ref{th:clt_reward}. This will be further
discussed in Section \ref{sec:clt_lipschitz}.

\begin{theorem}[Central limit theorem for costs]\label{th:clt_reward}
  Under  assumptions (A1,A2,A3,A4bis,A5), \\
  {\it (i)}- there exists constants $\beta $and $\gamma $ such that for all
  $x$:
  \begin{equation}
    \begin{array}{r}
      \displaystyle\limsup_{N\to\infty} \Pr( \sqrt{N}\Big|V^{*N}_{T}(\MM\toN_0,\CC\toN_0) -
      v^*_{T}(m_0,c_0)\Big| \geq x ) \\
      \leq  \Pr(\beta  \norminf{G_0}+ \gamma  \geq x);
    \end{array}
    \label{eq:clt_ii}
  \end{equation}
  {\it (ii)}-  there exist constants  $\beta', \gamma'>0$ such that for all $x$:
  \begin{equation}
    \begin{array}{r}
      \displaystyle\limsup_{N\to\infty} \Pr(\sqrt{N}\Big |V^{*N}_T(\MM\toN_0,\CC\toN_0) -
      V^N_{a_0^*\dots a^*_{T-1}}(\MM\toN_0,\CC\toN_0)\Big|\\ \geq x ) \leq
      \Pr(\beta' \norminf{G_0} + \gamma'\geq
      x);
    \end{array}
    \label{eq:clt_i}
  \end{equation}
  where  $\norminf{G'}=\sup_i |G'_i|$.
\end{theorem}

This theorem  is the main result of this section. 
The previous result (Theorem \ref{th:converg_opti}) says that  $ \limsup_{N\to\infty} V^{*N}_T(\MM\toN_0,\CC\toN_0)
=   \limsup_{N\to\infty} V^N_{a_0^*\dots a^*_{T-1}}(\MM\toN_0,\CC\toN_0) =  v^*_{t\dots
      T}(m_0,c_0)$.
This new theorem says  that both the gap   between the  cost under the
optimal policy and of the cost when using  the limit actions {\it
  (i)} or the gap between the latter cost and the optimal cost of the
limit system  {\it (ii)} are random variables  that decrease to 0 with speed
$\sqrt{N}$ and have  Gaussian laws.
Actually, a stronger result (using almost sure convergence instead of
convergence in law) will be shown in Corollary \ref{coro:rewardexp_clt}. A
direct consequence of this result is that there exists a constant
$\gamma''$ such that:  
\begin{equation}
  \E \Big[ \sqrt{N}|V^{*N}_{T}(\MM\toN_0,\CC\toN_0) -
  v^*_{T}(m_0,c_0)|\Big] \to \gamma''
\end{equation}

The rest of this section is devoted to the proof of this theorem.
A first step in the proof of Theorem  \ref{th:clt_reward} is a central limit theorem
for the states, which has an interest by its own.

\begin{theorem}[Mean field central limit theorem]
\label{th:clt_particles}
  Under assumption (A1,A2,A3,A4bis,A5-bis), if the actions taken by the
  controller are $a_0\dots a_{T-1}$, there exist Gaussian vectors of mean
  $0$, $G_1\dots G_{T-1}$ such that for every $t$:
  \begin{equation}
    \begin{array}{r}
      \sqrt{N}
      ((\MM\toN_0,\CC\toN_0) - (m_0,\ccc_0),
      \dots, (\MM\toN_t,\CC\toN_t) -
      (m_t,\ccc_t)
      ) \\
      \xrightarrow{\mathcal{L}} G_0,\dots, G_t.
    \end{array}
\label{eq:prop_clt}
  \end{equation}
  
  Moreover if $\Gamma_t$ is the covariance matrix of $G_t$, then: 
   \begin{equation} \Gamma_{t+1} = \left[
    \begin{array}{c|c}
      P_t&F_t\\
      \hline
      Q_t&H_t
    \end{array}
  \right]^{tr} \Gamma_t 
  \left[
    \begin{array}{c|c}
      P_t&F_t\\
      \hline
      Q_t&H_t
    \end{array}
  \right]
  + 
  \left[
    \begin{array}{c|c}
      D_t&0\\
      \hline
      0&0
    \end{array}
  \right] \label{eq:clt_relation_covaiance}
  \end{equation}
  where for all $1 \leq i,j \leq S$ and $1\leq k,\ell \leq d$:
  $(P_t)_{ij}{=}K_{ij}(a_t,\ccc_t)$, $(Q_t)_{kj}{=}\sum_{i=1}^S m_i \frac{\partial
    K_{ij}}{\partial c_k}(a_t,\ccc_t)$, $(F_t)_{ik}{=}\frac{\partial g_k}{\partial
    m_i}(m_{t+1},\ccc_t)$, $(H_t)_{k\ell} =\frac{\partial g_k}{\partial
    r_\ell}(m_t,\ccc_t)$,
  $(D_t)_{jj}=\sum_{i=1}^nm_i(P_t)_{ij}(1-(P_t)_{ij})$ and
  $(D_t)_{jk}=-\sum_{i=1}^nm_i(P_t)_{ij}(P_t)_{ik}$ ($j\neq k$).
\end{theorem}

\begin{proof}
  Let us assume that the Equation \eqref{eq:prop_clt} holds for some $t\geq
  0$. 
  
  As $\sqrt{N}((\MM\toN,\CC\toN)_t-(m,\ccc)_t)$ converges in law to $G_t$, there
  exists another probability space and random variables $\widetilde{\MM}\toN$ and
  $\widetilde{\CC}\toN$ with the same distribution as  $\MM\toN$ and $\CC\toN$ such that
  $\sqrt{N}((\widetilde{\MM}\toN,\widetilde{\CC}\toN)_t-(m,\ccc)_t)$ converges almost
  surely to $G_t$ \cite{durrett1991pta}. In the rest of the proof, by abuse of notation,
  we will write $\MM$ and $\CC$ instead of $\widetilde{\MM}$ and $\widetilde{\CC}$
  and then we assume that $\sqrt{N}((\MM\toN,\CC\toN)_t-(m,\ccc)_t) \toas G_t$.
  
  $G_t$ being  a Gaussian vector, there exists a vector of $S{+}d$
  independent Gaussian variables $U=(u_1,\dots,u_{S{+}d})^T$ and a matrix $X$
  of size $(S{+}d){\time}(S{+}d)$ such that $G_t = XU$.
  
  Let us call $P_t\toN\bydef K(a_t,\CC\toN_t)$.  According to lemma \ref{lem:sum_depent}
  there exists a Gaussian variable $H_t$ independent of $G_t$ and of
  covariance $D$ such that we can replace $\MM\toN_{t+1}$ (without changing
  $\MM_t$ and $\CC_t$) by a random variables
  $\widetilde{\MM}\toN_{t+1}$ with  the same laws such that:
  \begin{equation}
    \sqrt{N}(\widetilde{\MM}_{t+1}\toN-\MM_t\toN P_t\toN)\toas H_t.\label{eq:clt_eq1}
  \end{equation}
  In the following, by abuse of notation we write $\MM$ instead of 
  $\widetilde{\MM}$. Therefore we have
  \begin{equation*}
    \begin{array}{r}
    \sqrt{N}(\MM_{t+1}\toN{-}m_tP_t) = \sqrt{N}\Big( \MM_{t+1}{-}\MM\toN_t P\toN_t
    + m_t (P\toN_t{-}P_t) +\\ (\MM_t\toN {-} m_t)P_t +
    (\MM\toN_t{-}m_t)(P\toN_t{-}P_t)\Big)\\
    \toas H_t + m_t \displaystyle\lim_{N\to\infty} \sqrt{N}(P\toN_t{-}P_t)
    + \displaystyle\lim_{N\to\infty} \sqrt{N} (\MM_t\toN {-} m_t) P_t.
  \end{array}
  \end{equation*}

  By assumption, $\lim \sqrt{N}(\MM\toN_t-m_t)_i = (XU)_i$. Moreover,
  the first order Taylor expansion
  with respect to all component of $\CC$ gives  a.s. 
  \begin{eqnarray*}
    \lim_{N\to\infty}m_t\sqrt{N}(P\toN_t-P_t)_j &=& \sum_{i=1}^S m_{t_i}
    \sum_{k=1}^d \frac{\partial K_{ij}}{\partial c_{t_k}}(a_t, \ccc_t) (XU)_{S+k}\\
    &=& \sum_{k=1}^d Q_{kj} (XU)_{S+k}.
  \end{eqnarray*}

  Thus, the $j$th component of $\sqrt{N}(\MM_{t+1}\toN-m_tP_t)$ tends to
  \begin{equation}
    H_t + \sum_{k=1}^d Q_{kj} (XU)_{S+k} +
    \sum_{i=1}^S(XU)_iP_{ij}\label{eq:clt_sum_m_i} 
  \end{equation}

  Using similar ideas, we can prove that $\sqrt{N}(\CC\toN_{t_k}-c_{t_k})$ converges
  almost surely to $\sum_{i=0}^S\frac{\partial g_k}{\partial m_i} (XU)_i +
  \sum_{\ell =0}^d\frac{\partial g_k}{\partial c_{t_\ell}} (XU)_{S+\ell }$. Thus
  $\sqrt{N}((\MM\toN_{t+1},\CC\toN_{t+1})-(m_{t+1},\ccc_{t+1}))$ converges almost
  surely to a Gaussian vector. 

  Let us write the covariance matrix at time $t$ and time $t+1$ as two bloc
  matrices:  
  \[ \Gamma_t = \left[
    \begin{array}{c|c}
      \MM & O \\
      \hline
      O^T & \CC
    \end{array}
  \right]
  \mathrm{~and~} \Gamma_{t+1} = \left[
    \begin{array}{c|c}
      \MM' & O' \\
      \hline
      O'^T & \CC'
    \end{array}
  \right].
  \]
  For $1\leq j,j'\leq S$, $\MM'_{j,j'}$ is the expectation of
  \eqref{eq:clt_sum_m_i} taken in $j$ times \eqref{eq:clt_sum_m_i} taken in
  $j'$. Using the facts that $\E[(XU)_{S+k}(XU)_{S+k'}] = \CC_{kk'}$,
  $\E[(XU)_{S+k}(XU)_{i}] = O_{ik}$ and $\E[(XU)_{i}(XU)_{i'}]=\MM_{ii'}$,
  this leads to:
  \begin{equation*}
    \begin{array}{l}
      \begin{array}{r}
        \MM'_{j,j'} = \E[H_jH_j'] + \displaystyle\sum_{k,k'} Q_{kj}Q_{k'j'}\CC_{kk'} +
        \displaystyle\sum_{k,i'} Q_{kj}O_{i'k}P_{i'j'} \hfill \\+ \displaystyle\sum_{i,k'}
        Q_{k'j'}O_{ik'}P_{ij} + \displaystyle\sum_{i,i'} P_{ij}
        \MM_{ii'}P_{i'j}
      \end{array}\\
      = D_{jj'} {+} (Q^T\CC Q)_{jj'} {+} (Q^TO^TP)_{jj'} {+} (P^TOQ)_{jj'} {+}
      (P^T\MM P)_{jj'}.
    \end{array}
  \end{equation*}
  
  By similar computation, we can write similar equations for $O'$ and $\CC'$
  that lead to Equation \eqref{eq:clt_relation_covaiance}. 
\end{proof}

\begin{lem}\label{lem:sum_depent}
  Let $\MM\toN$ be a sequence of random measure on $\{1,\dots,S\}$ and $P\toN$
  a sequence of random stochastic matrices on $\{1,\dots,S\}$ such that
  $(\MM\toN,P\toN)\toas (m,p)$.  Let $(U_{ik})_{1\leq i \leq S,k\geq 1}$ be a
  collection of \emph{iid} random variables following the uniform 
  distribution on $[0;1]$ and independent of $P\toN$ and $\MM\toN$ and let us
  define $Y\toN$: for all $1\leq j\leq S$:
  \[ Y\toN_j \bydef \frac{1}{N} \sum_{i=1}^S \sum_{k=1}^{N\MM\toN_i}
  \mathbf{1}_{\sum_{l<k} P\toN_{il} < U_{ik} \leq \sum_{l\leq k}
    P\toN_{il}} \]
  then there exists a Gaussian vector $G$ independent of $\MM\toN$ and
  $P\toN$ and a random variable $Z\toN$ with the same law as $Y\toN$ such  
  that 
  \[\sqrt{N}(Z\toN - \MM\toN P\toN) \toas G.\]
  
  Moreover the covariance of the vector $G$ is the matrix $D$:
  \begin{equation}\label{eq:clt2_covariance}
    \left\{
      \begin{array}{llll}
        D_{jj} &=& \sum_i m_ip_{ij}(1-p_{ij})\\
        D_{jk} &=& -\sum_i m_ip_{ij}p_{ik} &(j\neq k).
      \end{array}
    \right.
  \end{equation}
\end{lem}

\begin{proof}
  As $(\MM\toN,P\toN)$ and $(U_{ik})_{1\leq i \leq S,k\geq 1}$ are independent, they can be
  viewed as functions on independent probability space $\Omega$ and
  $\Omega'$. For all
  $(\omega,\omega')\in\Omega{\time}\Omega'$, let
  $X_\omega\toN(\omega')\bydef \sqrt{N}
  (Y\toN(\omega,\omega')-\MM\toN(\omega) P\toN(\omega))$.
  
  By assumption, for almost all $\omega\in\Omega$, $
  (\MM\toN(\omega),P\toN(\omega))$ converges to $(m,p)$. A direct computation shows that,
  when $N$ grows, the characteristic function of $X^N_\omega$ converges to
  $\exp(-\frac{1}{2}\xi^T\sum_{i=1}^S m_iC_i\xi)$. Therefore for almost all
  $\omega$, $X\toN_\omega$ converges in law to $G$, 
  a Gaussian
  random variable on $\Omega'$.
  
  
  
  Therefore for almost all $\omega$, there exists a random variable
  $\widetilde{X}\toN_\omega$ with  the same law as $X\toN_\omega$ that
  converges $\omega'$-almost surely to $G(\omega')$. Let
  $Z\toN(\omega,\omega') \bydef \MM\toN(\omega)P\toN(\omega) +
  \frac{1}{N}\widetilde{X}\toN_\omega(\omega')$. By construction of
  $\widetilde{X}\toN_\omega$, for almost all $\omega$, $Z\toN(\omega,.)$
  has the same distribution as $Y\toN(\omega)$ and $\sqrt{N}(Z\toN-Y\toN
  P\toN)\xrightarrow{\omega,\omega'-a.s} G$. Thus there exists a
  function $\widetilde{Z}\toN(\omega,.)$ that has the same distribution as 
  $Y\toN(\omega)$ for all $\omega$ and that converges
  $(\omega,\omega')$-almost surely to $G$.  
\end{proof}

The first application of the mean field CLT is to show that it also works
for the cost. Let us assume that the controller takes actions $a_0\dots
a_{T-1}$ and let us introduce  the definition of $R\toN_{a_0\dots
  a_{T-1}}( \MM\toN_0,\CC\toN_0)=\sum_{t=1}^T(\re_t(\MM\toN_t,\CC\toN_t))$ and 
$r_{a_0\dots a_{T-1}}(m_0,\ccc_0)=\sum_{t=1}^T\re_t(m_t,\ccc_t)$. Lemma 
\ref{lem:converg_reward}, says that $R\toN_{a_0\dots a_{T-1}}( \MM\toN_0,\CC\toN_0) \toas
r_{a_0\dots a_{T-1}}(m_0,\ccc_0)$, the following results is more accurate:

\begin{coro}[\label{coro:reward_clt}Application of the CLT to reward]
  Under assumption (A1,A2,A3,A4-bis,A5-bis), if the controller takes
  the actions $a_0\dots a_{T-1}$ and if we call ${\bf D}r_t(m_t,\ccc_t)$ the
  differential of $\re_t(\MM,\CC)$ at the point $(m_t,\ccc_t)$, we have: 
  \begin{equation}
    \begin{array}{r}
      \sqrt{N}(R\toN_{a_0\dots a_{T-1}}( \MM\toN_0,\CC\toN_0)-r_{a_0\dots
        a_{T-1}}(m_0,\ccc_0))\\
      \tolaw\sum_{t=1}^T{\bf D}\re_t(m_t,\ccc_t)G_t.
    \end{array}
    \label{eq:reward_clt1}
  \end{equation}
\end{coro}

\begin{proof}
  Let $G_0\dots G_T$ be the Gaussian variables defined in the central limit  
  theorem. The proof of Theorem \ref{th:clt_particles} says that one
  can replace  $(\MM\toN_t,\CC\toN_t)$ by variables
  with  the same law such that the convergence is almost sure. Let $\omega$ be
  an event such that $\lim_N \sqrt{N}
  ((\MM\toN_t(\omega),\CC\toN_t(\omega))-(m,\ccc)_t))=G_t(\omega)$. For this
  event, we have $\lim_{N\to\infty} \sqrt{N}(\ccc_t(\MM\toN_t,\CC\toN_t)
  -\re_t(m_t,\ccc_t)) = {\bf D}\re_t(m_t,\ccc_t) G_t$ which leads to Equation
  \eqref{eq:reward_clt1} by using a Taylor expansion at order one.
\end{proof}

As the means of the Gaussian variables are 0, we have directly:

\begin{coro}\label{coro:rewardexp_clt}
  Under the same assumptions and if the convergence of the initial
  condition is almost sure ($( \MM\toN_0,\CC\toN_0)
  \toas(m_0,\ccc_0)$), one has:
  \begin{equation}
    \begin{array}{r}
      \sqrt{N}\Big|   V\toN_{a_0\dots
        a_{T-1}}( \MM\toN_0,\CC\toN_0)-v_{a_0\dots a_{T-1}}(m_0,\ccc_0)
      \Big| \\ \leq_{N\to\infty}  |{\bf D}\re_0(m_0,\ccc_t)G_0| \quad \mathrm{a.s.}
 \end{array}
    \label{eq:reward_clt2}
  \end{equation}
\end{coro}

\begin{proof}
  $v\toN_{a_0\dots a_{T-1}}(\MM\toN_0,\CC\toN_0)-v_{a_0\dots
    a_{T-1}}(m_0,\ccc_0) = \re(\MM\toN_0,\CC\toN_0)
-\re(m_0,\ccc_0)+\E_{\MM\toN_0,\CC\toN_0}[r\toN_{1\dots 
    T}(\MM\toN_1,\CC\toN_1) - \re_{1\dots T}(m_1,\ccc_1)]$. As
  $\sqrt{N}((\MM\toN_0, \CC\toN_0)-(m_0,\ccc_0))$ converges almost surely, the first
  part of the sum  can be upper bounded by
  $|{\bf D}\re_0(m_0,\ccc_0)G_0|$. 
  As for the second part of the sum, using the Berry-Esseen Theorem (Durrett
  2.4.d \cite{durrett1991pta}), one  can refine 
  Lemma \ref{lem:sum_depent} and show that the convergence is uniform.
Therefore  one can 
   switch the expectation and the limit, the second
  part of the sum  becomes
  $\E_{\MM\toN_0,\CC\toN_0}[\lim_{N\to\infty}\sqrt{N}(r\toN_{1\dots
    T}(\MM\toN_1,\CC\toN_1) - \re_{1\dots T}(m_1,\ccc_1))] =_{a.s} 0$
  which proves  Equation \eqref{eq:reward_clt2}.
\end{proof}

We are now ready for the  proof of  Theorem \ref{th:clt_reward}.

\begin{proof}[of theorem \ref{th:clt_reward}] 
  For a vector $G$, let us write
  $\normun{G}=\sum_i |G_i|$.
  Because of
   assumption (A4), there exists a compact set ${\mathcal B}$ such
   that for all $t$ from $0$ to $T$, 
   $\MM\toN_t,\CC\toN_t$ will remain in ${\mathcal B}$.

  Let us prove by induction on $t$ from $T$ to $0$ that there exist
  $\beta_t,\gamma_t\in\R^+$ such that if there exists a Gaussian variable $G_t$
  satisfying  $\sqrt{N}\big((\MM\toN_t,\CC\toN_t)-(m_t,\ccc_t)\big)\toas G_t$, then
  \begin{equation}
    \begin{array}{r}
  \limsup_{N\to\infty} \sqrt{N}\Big|V^{*N}_{t\dots T}(\MM\toN_t,\CC\toN_t)
  {-} v^*_{t\dots T}(m_t,\ccc_t)\Big|\\
  \leq  \beta_t \norminf{G_t}+\gamma_t.
  \end{array}
\label{eq:clt_preuve}
  \end{equation}
  
  For $t=T$, Corollary \ref{coro:rewardexp_clt}  can be used to transform
  Equation \eqref{eq:clt_preuve}  into 
  $\sqrt{N}|{\bf D}\re_T(m_T,\ccc_T)G_T| \leq
  \normun{{\bf D}\re_t(m_T,\ccc_T)} \norminf{G_T}$. Therefore,
  Inequality \eqref{eq:clt_preuve} is
  true if
  $\beta_T=\normun{{\bf D}\re_t(m_T,\ccc_T)}$ and  $\gamma_T=0$.
  
  Let us assume that \eqref{eq:clt_preuve} holds for some $t+1\leq T$ and that
  $\sqrt{N}\big((\MM\toN_{t},\CC\toN_{t})-(m_{t},\ccc_{t})\big)\toas G_{t}$.
  At time $t$, \eqref{eq:clt_preuve} can be upper bounded by: 
  \begin{equation*}
    \begin{array}{l}
      \sqrt{N}|\re_{t}(\MM\toN_{t},\CC\toN_{t}) -
      \re_{t}(m_{t},\ccc_{t})| \hfill \\ 
      \begin{array}{r}
      +  \sqrt{N}\Big|\sup_a\E_{\MM\toN_{t},\CC\toN_{t}}[V^{*N}_{t\dots
        T}(\Phi\toN_a(\MM\toN_{t},\CC\toN_{t}))] \\
      - \sup_a v^*_{t\dots T}
      (\Phi_a(m_{t},\ccc_{t}))\Big|.
    \end{array}
    \end{array}
  \end{equation*}

  The first part can be bounded by 
  $\normun{{\bf D}\re_t(m_t,\ccc_t)} \norminf{G_{t}}$. The rest of the  
  proof focuses in the second part of the sum. 
In the proof of Theorem \ref{th:clt_particles}, we showed that
  for all $a$ (up to the replacement  of  $\Phi\toN_a(\MM\toN_t,\CC\toN_t)$ by a
  random variable with  the same law), there exists a matrix $P_a$ and a
  Gaussian variable $G_a$ independent of $G_{t}$ such that
  $\sqrt{N}(((\MM\toN_{t},\CC\toN_{t}),(\MM\toN_{t+1},\CC\toN_{t+1}))-((m_{t},\ccc_{t}),(m_{t+1},\ccc_{t+1})))$
  converges almost surely to $(G_{t},P_aG_{t}+G_a)$.  Using the fact that
  $\sup_af(a)-\sup_ag(a) \leq \sup_a(f(a)-g(a))$, the expectation can be
  upper bounded by:
  \begin{equation*}
    \sup_a \sqrt{N} \E_{\MM\toN_{t},\CC\toN_{t}} \Big| V^{*N}_{t+1\dots
      T}(\Phi\toN_a(\MM\toN_{t},\CC\toN_{t})) - v^*_{t+1\dots T}(\Phi_a(m_{t},\ccc_{t}))\Big|.
  \end{equation*}
  
  Let us consider an arbitrary  action $a$. The Berry-Esseen Theorem shows that
  $\sqrt{N}((\MM\toN_{t+1},\CC\toN_{t+1})-(m_{t+1},\ccc_{t+1}))-P_aG_{t}$ converges uniformly to
  $G_a$, therefore we can switch the limit in $N$ and the expectation and
  by induction,  it can be upper bounded by
  $\E_G[\gamma_t\norminf{P_aG_{t}+G_a}+\beta_{t+1}] \leq \beta_{t+1} \norminf{P_aG_{t}} + \gamma_t +
  \beta_t\E[\norminf{G_a}]$. As $\A$ is compact and $(\MM\toN_{t+1},\CC\toN_{t+1})$
  remains in a compact set ${\mathcal B}$ (Equation \eqref{eq:unif_bounded}),
  $\sup_{a\in\A,(\MM,\CC)\in {\mathcal B}}\normun{P_a}<\infty$ and $\sup_{a\in\A,(\MM,\CC)\in
    {\mathcal B}}\E[\norminf{G_a}]<\infty$.  Thus to obtain an uniform bound on all
  $(\MM,\CC)$, taking $\beta_{t} \bydef \beta_{t+1}\sup_{\A,{\mathcal B}}\normun{P_a}$ and $\gamma_{t}\bydef
  \gamma_{t+1}+\beta_{t+1}\sup_{\A,{\mathcal B}}\E[\norminf{G_a}]$ satisfy \eqref{eq:clt_preuve}.

  Assumption (A4bis) says that at time $t=0$,
  $\sqrt{N}\big((\MM\toN_t,\CC\toN_t)-(m_t,\ccc_t)\big) \to  G_t$ holds in
  distribution. Using appropriate random variables $(\tilde{\MM}\toN_t,\tilde{\CC}\toN_t)$ with the
  same laws as  $(\MM\toN_t,\CC\toN_t)$ makes this convergence almost
  sure so that the induction above holds from $t=0$.
  This ends the proof for assertion {\it i} of the theorem. 
  
  As for assertion {\it ii},
  it comes from  the triangular inequality
  \begin{equation*}
    \begin{array}{r}
      \Big|V^{*N}_{T}(\MM\toN_0,\CC\toN_0) -
      V^{*N}_{a^*_0\cdots a^*_T}(\MM\toN_0,\CC\toN_0) \Big| \quad \quad \quad~ \\
      \begin{array}{l}
        \leq  \Big|V^{*N}_{T}(\MM\toN_0,\CC\toN_0) -
        v^*_{T}(m_0,\ccc_0)\Big| \\ ~\quad +\Big|  v^*_{T}(m_0,\ccc_0) - V^{*N}_{a^*_0\cdots
          a^*_T}(\MM\toN_0,\CC\toN_0) \Big|.
      \end{array}
    \end{array}
  \end{equation*}
  An upper  bound on the first term of the right side comes from
  assertion {\it i} and the second term can be bounded using Corollary
  \ref{coro:rewardexp_clt}.
  This ends the proof.
\end{proof}

\subsection{Infinite horizon discounted reward
  \label{sec:infinite-discouted}}

In this section, we prove  the first order results for 
infinite-horizon discounted Markov decision processes. As in the finite case, we
will show  that when $N$ grows large, the maximal expected discounted reward
converges to the one of the deterministic system and the optimal policy is
also asymptotically optimal. To do this  , we need the following
new assumptions: 

\begin{list}{(A\arabic{Lcount})}{\usecounter{Lcount}}
  \setcounter{Lcount}{5}
\item \label{ass:homo} \textbf{Homogeneity in time} -- The reward $\re_t$
  and the probability kernel $K_t$ do not depend on time: there exists
  $\re,K$ such that, for all $ \MM,\CC,a$ $\re_t(\MM,\CC) = \re(\MM,\CC)$
  and $K_t(a, \CC) = K(a,\CC).$
\item \label{ass:bounded} \textbf{Bounded reward} -- $\sup_{\MM,\CC} \re(\MM,\CC) \leq K <\infty.$
\end{list}

The homogeneity in time is clearly necessary as we are interested in
infinite-time behavior. Assuming that the cost is bounded might seems
strong but it is in fact very classical and holds in many situation, for
example when $\CC$ is bounded. The future reward are discounted according to a
discount factor $0\leq\delta <1$: if the policy is $\Pi$, the expected total
discounted reward of $\Pi$ is ($\delta$ is omitted in the notation):  
\[ V_\Pi\toN(\MM\toN_0,\CC\toN_0) \bydef  \E_\Pi\Big[\sum_{t=1}^\infty
\delta^{t-1}\re(\MM\toN_t,\CC\toN_t)\Big].\]
Notice that  Assumption (A7) implies that this sum remains
finite. The optimal total discounted reward $V^{*N}$ is the supremum
on all policies. 
For $T\in\N$, the
optimal discounted finite-time reward 
until $T$ is \[{V_{T}^*}\toN(\MM_0,\CC_0)  \bydef 
\sup_\Pi\E_\Pi\big[\sum_{t=1}^T  \delta^{t-1}\re(\MM_t,\CC_t)\big]. \]
As $\re$ is bounded, one can show that it
converges uniformly in $(\MM,\CC)$ to ${V^*}\toN$: 
\begin{equation}
  \lim_{T\to\infty} \sup_{\MM,\CC} \Big|{V_{T}^*}\toN(\MM,\CC) -
  {V^*}\toN(\MM,\CC) \Big| = 0.\label{eq:discounted_fini} 
\end{equation}
Equation \eqref{eq:discounted_fini} is the key of the following analysis. Using
this fact, we can prove the convergence when $N$ grows large for fixed $T$
and then let $T$ go to infinity. Therefore with a very few changes in
the proofs of Section \ref{sec:finite-horizon}, we have the following result:
\begin{theorem}[Optimal discounted case]
  Under assumptions (A1,A2,A3,A4,A6,A7), as $N$ grows large, the optimal
  discounted reward of the stochastic system converges to the optimal
  discounted reward of the deterministic system:
  \begin{equation*}
    \lim_{N\to\infty} {V^*}\toN(\MM\toN,\CC\toN) =_{a.s} v^*(m,c),
  \end{equation*}
 where $v^*(m,c)$ satisfies the Bellman equation for the deterministic
 system:
  \[  v^*(m,c) = \re(m,c) + \delta \sup_{a\in\A} \Big\{v^*(\Phi_a(m,c)) \Big\}.\]
\end{theorem}

\subsubsection{Problems for other infinite horizon criteria}

Again, the discounted problem is very similar to the finite case because
the total reward mostly depends on the rewards during a finite amount of
time. As for other other infinite-horizon criteria such as average reward
or its variants, the average reward is (if it exists) $\lim_{T\to\infty}
\frac{1}{T} \E_\Pi \sum_{t=1}^T c(\MM_t,\CC_t)$
. 

This raises the problem of the exchange of the limits $N\to\infty$ and
$T\to\infty$. Consider a case without control with two states
$\SS{=}\{0;1\}$ and $\CC_t$ is the mean number of particles in state $1$
($\CC_t=(\MM_t)_1$) and with a function $f{:}[0;1]{\to}[0;1]$ such that the
transition kernel $K$ is $K_{i1}(\CC)=f(\CC)$ for $i\in\SS$. If
$\MM^N_0(0)\toas m_0$ then for any fixed $t$, $\MM^N_t$ converges to
$f(f(\dots f(m_0)\dots))$. Using
techniques that can be found in \cite{borkarStochasticApprox}, one can
prove that as $N$ grows large, $\lim_{t\to\infty}\MM\toN_t$ might converges
to almost any subset of $L{\subset}[0;1]$ such that $L=f(L)$. However, in
general $\lim_{t\to\infty}\lim_{N\to\infty}\MM\toN_t \neq
\lim_{N\to\infty}\lim_{t\to\infty}\MM\toN_t$. For example if
$f(x)=x$, the deterministic system is constant while the stochastic system
converges almost surely to a random variable (as a bounded Martingale) that
takes values in $\{0;1\}$. 

Similar difficulties arise for the
central limit theorem in the discounted case: the convergence depends on the behavior of the
system when $T$ tends to infinity. 

\section{Application to a brokering problem} 
\label{sec:example}

To illustrate the usefulness of our framework, let us consider the
following model of a brokering problem in computational grids.
There are $A$ application sources that send 
tasks  into a grid  system and a central broker routes all
theses tasks  into $d$ clusters 
(seen as multi-queues) and tries to minimize the total waiting time of the tasks.
A similar queuing model of a grid broker was used in \cite{Mitrani,ArticleBerten.BG_PC07,InProceedingsBerten.BG_gbfbau_07}.

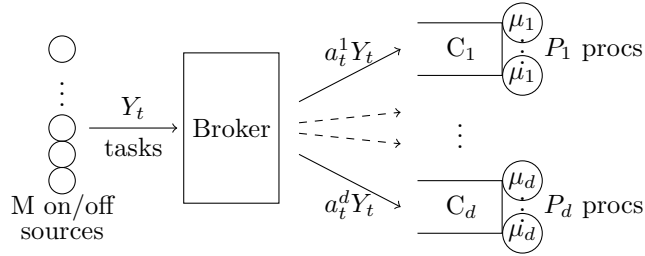
\begin{figure}[ht]
  \centering
  \begin{tikzpicture}[xscale=0.9, yscale=0.7]
    \foreach \x in {1,2,3,6} \node[circle,draw,minimum size=10] (\x) at (3.5,.5*\x-1.5) {};
    \node at (3.5,.8) {$\vdots$};
    \node at (3.5,-1.5) {M on/off};
    \node at (3.5,-2) {sources};
    \draw (3.9,0) edge[->] node[above] {$Y_t$} node[below] {tasks} (5.2,0);
    \draw (7,.5) edge[->] node[above] {$a_t^1Y_t$} (8.5,1.5);
    \draw (7,-.5) edge[->] node[below] {$a_t^dY_t$} (8.5,-1.5);
    \draw (7,.1) edge[->,dashed] (8.5,.3);
    \draw (7,-.1) edge[->,dashed] (8.5,-.3);
    \node[rectangle, draw, minimum height=2cm] at (6.,0) {Broker};
    \foreach \x/\y in {1.5/1,-1.5/d}
    \draw (8.75, \x+.5) -- (10, \x+.5) -- (10., \x-.5) -- (8.75, \x-.5)
    node at (10.3,\x+.1) {$\vdots$}
    node[circle, draw,minimum size=15, inner sep=0pt] at (10.3,\x+.5) {$\mu_\y$}
    node[circle, draw,minimum size=15, inner sep=0pt] at (10.3,\x-.5) {$\mu_\y$}
    node at (11.35,\x) {$P_\y$ procs}
    node at (9.4,\x) {$\CC_\y$};
    \node at (9.38,0) {$\vdots$};
  \end{tikzpicture}
  \caption{The routing system}
  \label{fig:routing}
\end{figure}

Here, time is discrete and the $A$ sources follow a discrete on/off model:
for each source $j\in\{1\dots A\}$, let $(Y_t^j)\bydef 1$ if the source is \textrm{on} (\emph{i.e.} it
sends a tasks between $t$ and $t+1$) and $0$ if it is \textrm{off}.  The
total number of packets sent between $t$ and $t+1$ is $Y_t\bydef \sum_j
Y_t^j$. Each queue $i\in\{1\dots d\}$ is composed of $P_i$ processors, and
all of them work at speed $\mu_i$ when available. Each processor $j\in\{1\dots P_i\}$ of the queue
$i$ can be either \emph{available}  ( in that case we set 
$X_t^{ij}\bydef 1$ )  or \emph{broken} (in that case $X_t^{ij}\bydef
0$). The total number of processors available in the queue $i$ between $t$ and
$t+1$ is $X_t^i\bydef \sum_j X_t^{ij}$ and we define $B_t^i$ to be the
total number of tasks waiting  in the queue $i$ at time $t$. At each time slot
$t$, the broker (or controller) allocates  the $Y_t$ tasks to the $d$
queues: it chooses an action $a_t\in\mathcal{P}(\{1\dots Y_t\}^d)$ and
routes each $Y_t$ packets in queue $i$ with probability $a^i_t$. The system
is represented figure \ref{fig:routing}.  The number of tasks in the queue
$i$ (buffer size) evolves according to the following relation: 
\begin{equation}
  B_{t+1}^i = \Big(B_{t}^i - \mu_i X_t^i + a_t^iY_t \Big)^+.\label{eq:routing_ex:queue}
\end{equation}

The cost that we want to minimize is the sum of the   waiting times of
the tasks. Between $t$ and $t+1$, there are $\sum_iB^t_i$ tasks
waiting in the queue, therefore the cost at time $t$ is
$r_t(B)\bydef\sum_i{B^i_t}$.  As we consider a  finite horizon, we should
decide a cost for the remaining tasks in the queue. In our simulations,
we choose $r_T(B)\bydef\sum_i{B^i_T}$. 

This problem can be viewed as a multidimensional restless bandit
problem where computing the optimal policy for the broker is
known to be a hard problem \cite{Whittle}. Here, indexability may help to compute
near optimal policies by solving one MDP for each queue \cite{Whittle,Weber}. However the
complexity remains high when the number of processors in  all the queues and the number of
sources are large.

\subsection{Mean field limit}

This  system can be modeled using the  framework of particles evolving in a common
environment. 

  \begin{list}{\textbullet}%
    {\setlength{\leftmargin}{0.5\leftmargin}}
\item There are $N  \bydef
A+\sum_{i=1}^dP_i$ ``particles''. Each particle can either be a source
(of type $s$) or a server (belonging to one of the queues, $q_1\cdots q_d$),
and can either be ``on'' or ``off''. Therefore, the possible states
of one  
particle is an element of $\SS=\big\{ (x,e) | x \in \{s, q_1, \cdots ,q_d\} , e \in
\{\mathrm{on},\mathrm{off}\} \big\}$. the population mix $M$ is the
proportion of sources in state on and the proportion of servers in
state on, for each queue.

\item The action of the controller  are the routing choices of the broker:
$a_t^d$ is the probability that a task is  sent to queue $d$ at time
$t$.

\item The environment of the system depends on the vector $B_t = (B_{t_1}\dots B{t_d})$,
giving the number of tasks in queues $q_1,\dots q_d$ at time $t$.
The time evolution of the i-th component  is 
\begin{equation*}
B_{t+1_i}  = g_i(B_t, \MM\toN_{t+1}, a_t) \bydef  \Big(B_{t_i}- \mu_i X_t^i + a_t^iY_t \Big)^+.
\end{equation*}
 The shared environment is represented by the context
 $C^N_t \bydef  (\frac{B_{t_1}}{N}\dots\frac{B_{t_d}}{N})$.

\item Here, the transition kernel can be time dependent but is independent
  of $a$ and $C$.  The probability of a particle to go from a state
  $(x,e)\in\SS$ to $(y,f)\in\SS$ is $0$ if $x \not= y$ (a source cannot
  become a server and vice-versa). If $x = y$ then $K_{(x,
    \mathrm{on}),(x,\mathrm{off})}(a,C)(t)$ as well as
  $K_{(x,\mathrm{off}),(x,\mathrm{on})}(a,C)(t)$
  are arbitrary
  probabilities.

\end{list}

Here is how a system of size $N$ is defined.
A preliminary number of sources  $A_0$ as well as a preliminary
number  $P_i$ of servers per queue is given, totaling in $N_0$ particles.
For any  $N$, a   system with  $N$ particles is composed  of $\lfloor A_0N/N_0 \rfloor$ (resp. $\lfloor
P_iN/N_0 \rfloor$) particles that are sources (resp. servers  in queue
$i$). The remaining particles (to reach a total of  $N$) are allocated
randomly with a
probability proportional to the fractional part of $A/N_0$ and $P_iN/N_0$
so that the mean number of particles that are sources is $A/N_0$ and
the   mean number of particles that are servers  in queue
$i$ is $P_iN/N_0$. Then, each
of these particles  changes state  over time  according to 
the probabilities $K_{u,v}(a,C)(t)$. At time $t=0$, 
a particle is in state ``on'' with probability one half.

It should be clear that this system satisfies
Assumptions (A1) to (A4) and therefore one can apply the convergence
theorem \ref{th:converg_opti} to this system that shows that if using 
the policies $a^*$ or $\Pi^*$, when $N$ goes to infinity the 
system  converges  to a deterministic system with optimal cost. 
An explicit computation of the policies $a^*$ and $\Pi^*$ is possible
here and is postponed to Section \ref{sec-compuOPt}.

\subsection{CLT applicability}
\label{sec:clt_lipschitz}

As for the central limit theorem, 
Assumption (A4-bis) on the convergence of the initial condition to a
Gaussian variable is true since the random part of the initial state is
bounded by $ \frac{N_0}{N}$ and $\sqrt{N} \frac{N_0}{N} $ goes to 0 as $N$
grows.  Unfortunately Assumption (A5) does not hold since the
function $g$ is not
differentiable when $\CC_{t}^i{-}\mu_i X_t^i{+}a_t^iY_t=0$.  However, as
mentioned  in the beginning of section \ref{sec:clt} the
differentiability condition in
Assumption (A5) can be replaced by a Lipschitz continuity condition. Let us
consider Assumption (A5-ter):

\begin{list}{}{}
\item[(A5-ter)] \textbf{Continuous Lipschitz} -- For all $t$ and all
  $i,j\in\SS$, all functions $g$, $K_{ij}$ and $\re_t$ are Lipschitz
  continuous on all compact sets  of their domain.
\end{list}

This assumption is weaker than (A5) since, if a function is $C^1$, it is
Lipschitz on every compact set  (with Lipschitz constant $\sup||f'||$). In
the example, function $g$ has a 
right-derivative and a left-derivative at all points and therefore satisfies
(A5-ter). The central limit theorem \ref{th:clt_reward} should apply here
as well: 

\begin{theorem}
  \label{th:clt_lipschitz} Theorem \ref{th:clt_reward} still holds
  when replacing (A5) by (A5-ter).
\end{theorem}

\begin{proof}[(Sketch of the proof)] The proof is very similar to the one of
  \ref{th:clt_reward} and we just sketch  the main differences.
  
  As seen at the end of section \ref{sec:assum}, all variables are almost
  surely bounded. By assumption (A5-ter), all functions are Lipschitz, thus
  let $L_g,L_K,L_{\re_t}$ be the Lipschitz constants on the compact space
  $\mathcal{B}$ (see Equation \eqref{eq:unif_bounded}) for $g,K$ and $\re_t$ respectively
  and $L=\max\{L_g,L_K,L_{\re_t}\}$. The main idea is to replace all
  equalities in the proof of all CLT theorems by inequalities. For
  instance, in Theorem \ref{th:clt_particles}, Equation \eqref{eq:prop_clt}
  is replaced by the following statement: for all $x_1\dots x_t\in\R^t$,
  \begin{equation}
    \begin{array}{r}
      \limsup_{N}\Pr\big(\sqrt{N}(\norminf{(\MM\toN_0,\CC\toN_0) -
        (m_0,\ccc_0)},\dots, \quad\quad ~ \\\norminf{(\MM\toN_t,\CC\toN_t) -
        (m_t,\ccc_t)}\big) \geq (x_1\dots x_t)\big) \quad~\\
      \leq  \Pr(
      (\norminf{G_0},\dots, \norminf{G_t}) \leq (x_1\dots x_t))
    \end{array}
  \end{equation}
  where the variables $G_t$ have covariance $\Gamma_t = L^2
  \Gamma_{t-1}+D_{t-1}$. The other steps in the proof can  be changed in
  almost the same way.  
Formula \eqref{eq:reward_clt1} in Corollary \ref{coro:reward_clt} is  replaced by 
  \begin{equation}
    \begin{array}{r}
    \sqrt{N} | R\toN_{a_0\dots a_{T-1}}(
    \MM\toN_0,\CC\toN_0)-r_{a_0\dots a_{T-1}}(m_0,\ccc_0)| \\ \leq_{st}
    \sum_{t=0}^T  L \norminf{G_t}
  \end{array}
    \label{eq:reward_clt3}
  \end{equation}
and Formula \eqref{eq:reward_clt2} of Corollary \ref{coro:rewardexp_clt} by 
  \begin{equation}
    \begin{array}{r}
      \sqrt{N}\Big|  V\toN_{a_0\dots a_{T-1}}(
    \MM\toN_0,\CC\toN_0)-v_{a_0\dots a_{T-1}}(m_0,\ccc_0)     \Big| \\
   \leq   \alpha \norminf{G_0} + \delta , \quad \mathrm{a.s.}
 \end{array}
    \label{eq:reward_clt4}
  \end{equation}
where $\alpha$ and $\delta  $ are constants depending  on $L$.
\end{proof}

\subsection{Optimal policy for the deterministic limit}
\label{sec-compuOPt}

As the evolution of the sources and of the processors does  not depend  on the
environment, for all $i$, $t$, the quantities $\mu_iX_t^i$ and $Y_t$
converge 
almost surely to deterministic values that we call $x^i_t$ and $y_t$. If 
$y^i_t$ is the number of packets distributed to the $i$th queue at time
$t$,  $c^i_{t+1} = (c^i_t+y^i_t-x^i_t)^+$. The deterministic optimization problem is
to compute 
\begin{equation}
  \label{eq:routing_ex:deterministic}
  \min_{y^1_1\dots y^d_T}\{ \sum_{t=1}^T\sum_{i=1}^d c^i_t \mathrm{~with~}
  \begin{array}{l}
    c^i_{t+1} = (c^i_t+y^i_t-x^i_t)^+\\
    \sum_i y^i_t = y_t
\end{array}
\}.
\end{equation}

Let us call $w^i_t$ the work done by the queue $i$ at time $t$:
$w^i_t=c^i_{t}-c^i_{t-1}+y^i_{t-1}$.  The sum of the size of the queues at
time $t$ does not depend on with queue did the job but only on the quantity
of work done:
\[ \sum_{i=1}^d c^i_t = \sum_{i=1}^d c^i_0 - \sum_{u\leq t,i} w^i_t \]
Therefore to minimize the total cost, we have to maximize the total work
done by the queues. Using this fact, the optimal strategy can be computed
by iteration of a greedy algorithm.

\begin{figure}[ht]
  \centering
  \definecolor{light-gray}{gray}{0.85}
  \newcommand\off{\cellcolor{light-gray}}
    \begin{tabular}{|c|c|c|c|c|c|c|c|c|}
      \hline
      Time $t$ & 0 & 1 & 2 & 3 & 4 & 5 & 6\\
      \hline \hline
      $y_t$ (tasks)
      & 8  & 1  & 0  & 1   & 7  &  6 & 6\\
      \hline \hline 
      \multirow{4}{*}{Queue 1}
      &  X & $\tau_0$ & $\tau_0$ & $\tau_3$   & $\tau_4$ & $\tau_4$ & $\tau_6$ \\ \cline{2-8}
      &  X & $\tau_0$ & $\tau_0$ &     & $\tau_4$ & $\tau_5$ &\off\\ \cline{2-8}
      &  $\tau_0$&\off&\off&\off & $\tau_4$ &\off&\off\\ \cline{2-8}
      &\off&\off&\off&\off & $\tau_4$ &\off&\off\\ \cline{2-8}
      \hline
      \multirow{2}{*}{Queue 2}
      &\off&\off & X & X   &\off& $\tau_5$ & $\tau_6$ \\ \cline{2-8}
      &\off&\off&\off&     &\off& $\tau_5$ &\off\\ \cline{2-8}
      \hline
      \multirow{3}{*}{Queue 3}
      &  X & $\tau_0$ & $\tau_1$ &     & $\tau_4$ & $\tau_5$ & $\tau_6$ \\ \cline{2-8}
      &  X & $\tau_0$ &    &     & $\tau_4$ & $\tau_5$ & $\tau_6$  \\ \cline{2-8}
      &\off& $\tau_0$ &    &     &\off& $\tau_5$ &\off\\ \cline{2-8}
      \hline \hline
      \multirow{3}{*}{Optimal allocation}
      &5&.&.&1&5&1&1+2\\ \cline{2-8}
      &.&.&.&.&.&2&1\\ \cline{2-8}
      &3&1&.&.&2&3&2\\ \cline{2-8}
      \hline
    \end{tabular}
  \caption{This figure presents an example of an execution of the
    algorithm. We consider a case with $3$ queues. At $t=0$ (resp. $1,...,6$)
    there are $8$ (resp. $1,0,1,7,6,6$) packets arriving in the system. Each
    processor has speed $1$ and the  processors in state
    ``off'' are represented by grey cells (for example, at time $0$, there
    are respectively $3,0$ and $2$ processors available in queue $1,2$ and
    $3$). All queues start at time 
    $0$ with $2$ packets. The top part of the table shows at which time a
    packet will be processed while the bottom part shows the corresponding
    optimal allocation (\texttt{X} represent tasks present  in the queues
    before $t=0$; A label $\tau_i$ in a slot of queue $j$ at time $t$ represents
    one  task arriving at time $i$ allocated to queue $j$ that will be processed
    at time $t$. The number of slots with label $\tau_i$ should be equal to $y_i$;
    At the end, $2$ packets cannot be allocated in empty
    slots. They are  routed arbitrarily  (in queue $1$)).
}
    \label{fig:compute_determi_policy}
\end{figure}
The principle  of the algorithm is the following.
\begin{enumerate}
\item The processors in all queues,  which are ``on'' at time $t$  with a
speed $\mu$ are seen as slots of size $\mu$.

\item
At each  time $t$,  $y_t$ units of tasks have to be allocated. This is done
in a greedy fashion by filling up the empty slots starting from time
$t$.
Once all slots at time $t$ are full, slots at time $t+1$ are
considered and are filled up with the remaining volume of tasks, and
so forth up to time $T$.

\item
The remaining tasks  that do not fit in the slots before $T$ are
allocated in an arbitrary fashion.
\end{enumerate}

See figure \ref{fig:compute_determi_policy} for an illustration of the
execution of the algorithm on an example.  It should be clear that the
algorithm is linear in the number of slots $nk$ and that this algorithm
computes an optimal allocation.

\subsection{Numerical example}

\begin{figure}[ht]
  \centering
  
  \begin{tabular}{cc}
    \begin{tabular}{c}
      \rotatebox{90}{Cost}
    \end{tabular}
    &
    \begin{tabular}{c}
      \hspace{-1cm}
      \includegraphics[angle=-90,width=.9\linewidth]{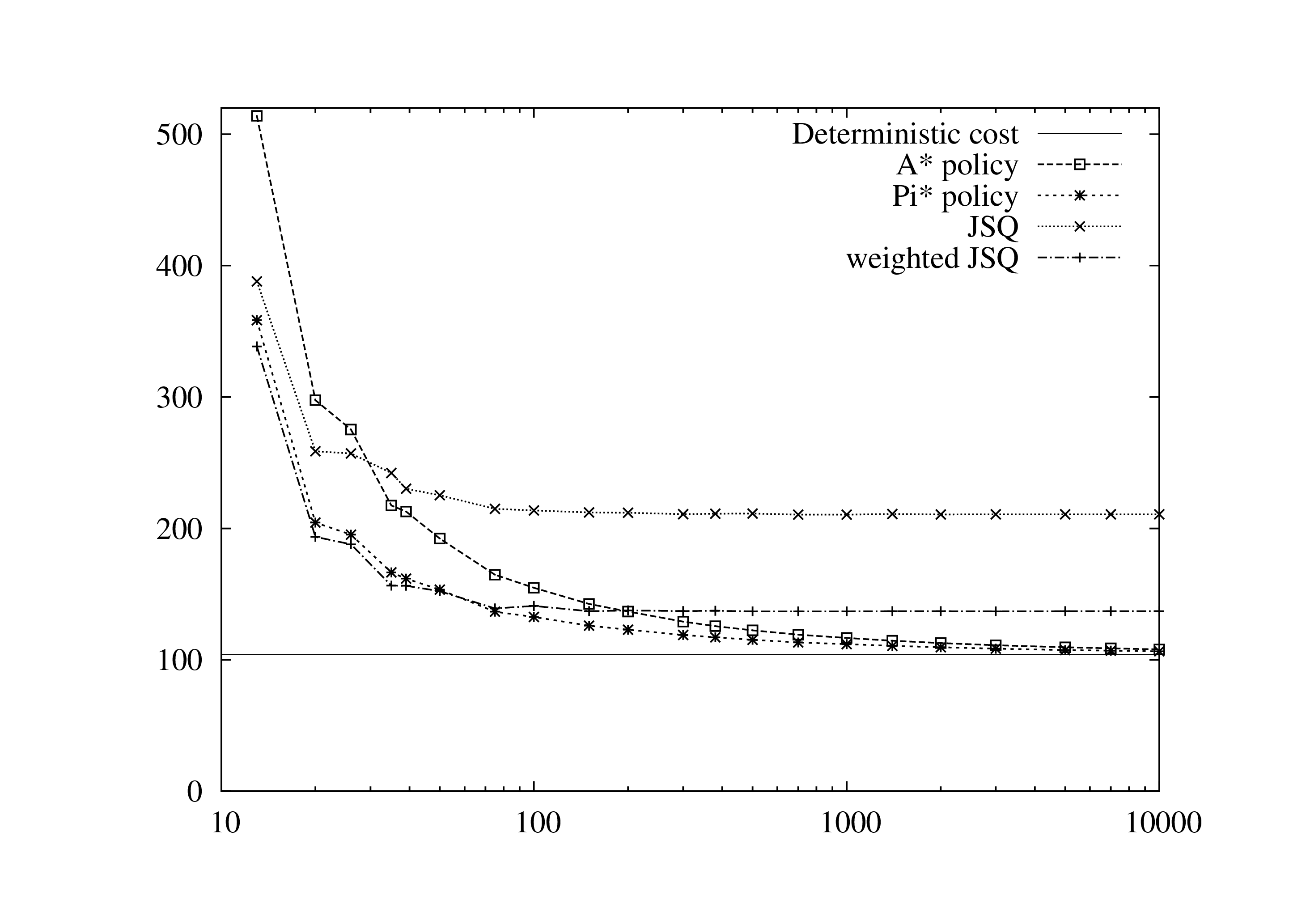}\\
      Size of the system: $N$ 
    \end{tabular}
  \end{tabular}
  \caption{Expected cost of the policies $a^*$, $\Pi^*$, JSQ and W-JSQ for
    different values of $N$. \label{fig:conv}}
\end{figure}

\begin{figure}[ht]
  \centering

  \begin{tabular}{cc}
    \begin{tabular}{c}
      \rotatebox{90}{$\sqrt{N}(V^N_X-v^*)$}
    \end{tabular}
    &
    \begin{tabular}{c}
      \hspace{-1cm}
      \includegraphics[angle=-90,width=.9\linewidth]{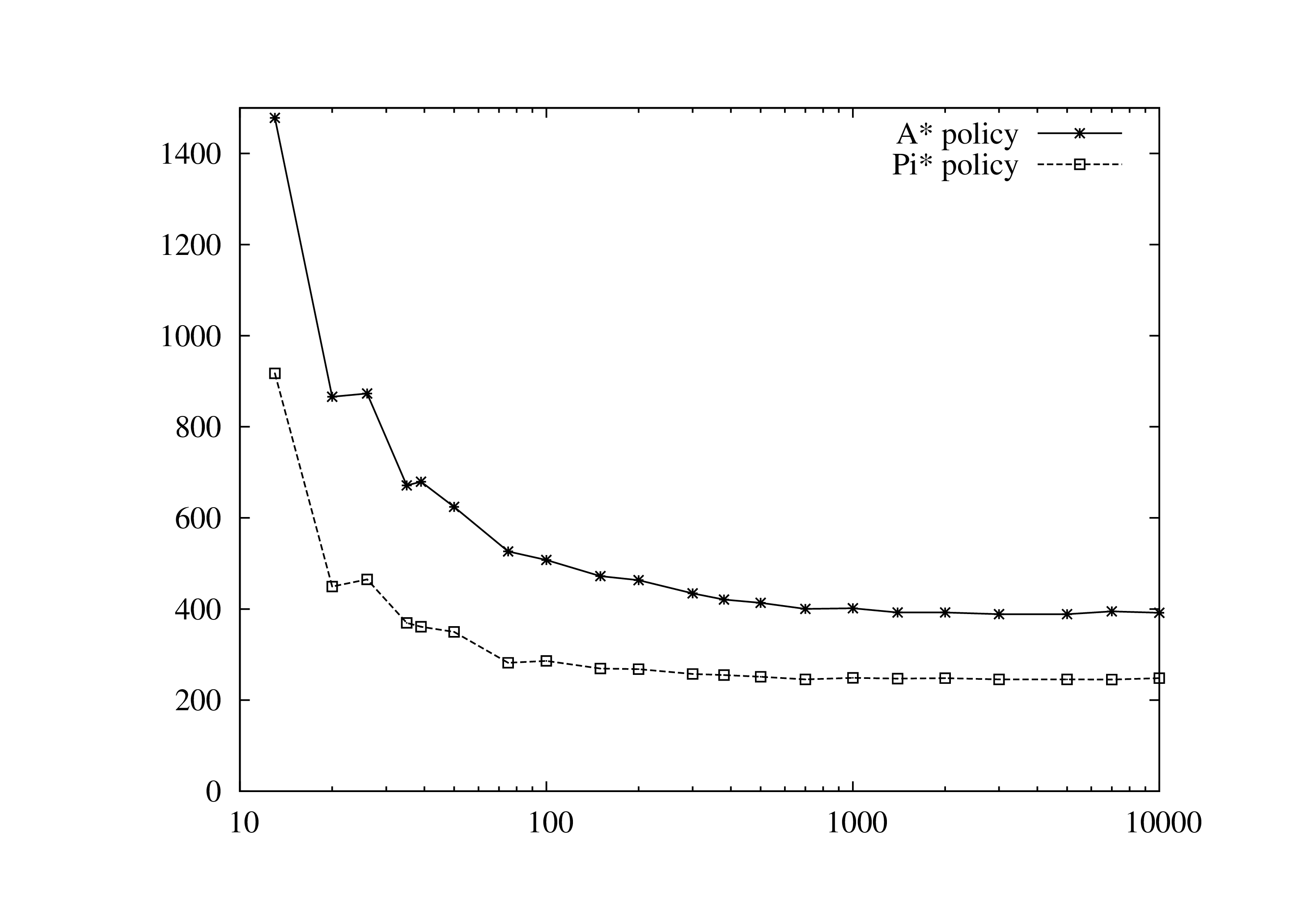}\\
      Size of the system: $N$
    \end{tabular}
  \end{tabular}
 \label{fig:convergence_simulation}

  \caption{Speed of convergence  of the policies $X = a^*$ or $\Pi^*$ for
    different values of $N$. \label{fig:conv_clt} }
\end{figure}

We consider a simple instance of the resource allocation problem with
$5$ queues. Initially, 
they  have respectively $1,2,2,3$ and $3$
processors running at speed $.5,.1,.2,.3$ and $.4$ respectively. There are
$3$ initial  sources. The transition matrices are time dependent and are
chosen randomly before the execution of the algorithm -- that is they are
known for the computation of the optimal policy and are the same for all
experiments.  We ran some simulations to compute the expected cost of
different policies for various sizes of the system. We compare different
policies:
\begin{enumerate}
\item Deterministic policy $a^*$ -- to obtain this curve, the optimal
  actions $a^*_0\dots a^*_{T-1}$ that the controller must take for the
  deterministic system have been computed.  At time $t$, action $a^*_t$ is
  used regardless of the currently state, and the cost up to time $T$ is
  displayed.
\item Limit policy $\Pi^*$ -- here,  the optimal policy $\Pi^*$ for
  the deterministic case was first computed. When the stochastic
  system is in state  $(\MM^N_t,\CC^N_t)$ at
  time $t$, we apply the action $\Pi^*_t(\MM^N_t,\CC^N_t)$ and the
  corresponding cost up to time $T$ is reported.
\item Join the Shortest Queue (JSQ) and Weighted Join the Shortest Queue (W-JSQ) --
  for JSQ, each packet is routed (deterministically) in the shortest
  queue. In W-JSQ, a packet is routed in the queue whose weighted  queue size
  $B_i/(\mu_iX_i)$ is the smallest.
\end{enumerate}
The results are reported in Figures \ref{fig:conv} and \ref{fig:conv_clt}.

A series of several simulations for with different values of $N$ was run.
The reported values in the figures are the mean values of the waiting time
over 10000 simulations for small values of $N$ and around 200 simulations
for big values of $N$. Over the whole range for $N$, the 95\% confidence
interval is less than 0.1\% for the expected cost -- figure \ref{fig:conv}
-- and less than 5\% for the central limit theorem -- figure
\ref{fig:conv_clt}.

Figure \ref{fig:conv} shows the average waiting time  of the
stochastic system when we apply the different policies. The horizontal
line represents the optimal cost of the deterministic system
$v^*(m_0,c_0)$ which is probably less than $V^{*N}(\MM_0,\CC_0)$. This 
figure illustrates Theorem \ref{th:converg_opti}: if we apply $a^*$ or 
$\Pi^*$, the cost converges to $v^*(m_0,c_0)$. 

In Figure \ref{fig:conv}, one can see that for low values of $N$, all the curves are
not smooth. This behavior comes from the fact  that 
when $N$ is not very large with respect to $N_0$, there are at least $\lfloor\frac{N}{N_0}A\rfloor$ 
(resp. $\lfloor \frac{N}{N_0}P_i\rfloor$) particles that are sources
(resp. processors in queue $i$) and the remaining particles are distributed
randomly. The random choice of the remaining states are chosen so that
$\E[A\toN]= \frac{N}{N_0}A$, but the difference $A^N  -  {N}{N_0}A$
may be large. Therefore,  for some $N$ the load of the system is much
higher than the average load, leading to larger costs.
As $N$ grows, the
proportion of remaining particles decreases and the phenomena becomes
negligible. 

A second feature that shows in  Figure \ref{fig:conv},
is the fact that on all curves, the expected waiting times are decreasing when
$N$ grows.
This behavior is certainly related to Ross conjecture \cite{Rolski} that
says that for a given load,  the  average queue length
decreases when the arrival and service processes are more deterministic.

Finally, the most important information on this figure is the fact that the
optimal deterministic policy and the optimal deterministic actions perform
better than JSQ and weighted JSQ as soon as the total number of elements in
the system is over 200 and 50 respectively.  The performance of the
deterministic policy $a^*$ is quite far from W-JSQ and JSQ for small values
of $N$, and it rapidly becomes better than JSQ ($N\geq 30$) and W-JSQ
($N\geq 200$). Meanwhile the behavior of $\Pi^*$ is uniformly good even for
small values of $N$.

The figure \ref{fig:conv_clt} illustrates Theorem \ref{th:clt_reward} 
which says that the speed of convergence towards the limit is of order
$\sqrt{N}$. On the $y$-axis, $\sqrt{N}$ times the
average cost of the system minus the optimal deterministic cost is
plotted.  One can see
that the gap between the expected cost of the policy $\Pi^*$
(resp. $a^*$) and the deterministic cost $v^*(m_0,c_0)$ is about
$250/ \sqrt{N}$ (resp. $400 / \sqrt{N}$) when $N$ is large.This should
be an upper bound on the constant $\delta$ defined in Equation \eqref{eq:reward_clt4}.

Besides comparing $a^*$ and $\Pi^*$ to other heuristics, it would be
interesting to compare it to the optimal policy of the stochastic
system, whose cost is 
$V^{*N}(\MM,\CC)$. One way to compute this optimum would be by using Equation
\eqref{eq:v_*_1...t}. However to do so, one needs  to solve  it for all
possible values of $\MM$ and $\CC$. In this example, $\CC$ can be as
large as  the length of
the five queues and each particle's state can vary in \{on,off\}. Therefore
even with $N=10$ and if we only compute the cost for queues of size less
than 10, this leads to $2^N10^5 \approx 10^8$ states which is hard to
handle even with powerful computers.

\section{Computational issues}

Throughout  the paper, we have shown that if the controller uses  the optimal
policy $\Pi^*$ of the deterministic limit of the finite real system, the
expected cost will be close to the optimal one (Theorem
\ref{th:converg_opti}). Moreover, Theorem \ref{th:clt_reward} gives
a bound on the error that we make. However to apply these results in
practice, a question remains: how difficult is it to compute the
optimal limit policy? 

The first answer comes straight from the example. In many cases, even if
the stochastic system is extremely hard to solve, the deterministic limit
is often much simpler. The best case of course is, as in the example of
section \ref{sec:example}, when one can compute the optimal policy. If one
can not compute it, there might also exist approximation policies with
bounded error (see \cite{hochbaum1996aan} for a review on the
subject). Imagine that a 2-approximation algorithm exists for the
deterministic system, then, Theorem \ref{th:converg_opti} proves that for
all $\varepsilon$, this algorithm will be a
$(2{+}\varepsilon)$-approximation for the stochastic system if $N$ is large
enough. Finally, heuristics for the deterministic system can also be
applied to the stochastic version of the system.

If none of this works properly, one can also compute  the optimal
deterministic policy by  ``brute-force'' computations 
using Equation \eqref{eq:v_*_1...t}: $v^*_{t\dots T}(m,c) =
\re_t(m,c)+\sup_{a} v^*_{t+1\dots
  T}(\Phi_a(m,c))$. In that case,  an approximation
of the optimal policy is obtained by discretizing the state space
and by solving the equation backward (from $t=T$ to $t=0$), to obtain 
the optimal policy for all states.  The
brute force approach 
can also be applied directly on the stochastic equation using \eqref{eq:V_*_1...t}:
$V^{*N}_{t\dots T}(\MM,\CC)=\re_t(\MM,\CC)+\sup_{a\in\A} \E_{\MM,\CC}
\Big[V^{*N}_{t+1 \dots T}\big(\Phi\toN_a(\MM,\CC)\big)\Big]$. However,
solving  the deterministic system has three  key advantages.  
The first one is that the size of the discretized deterministic system
may have nothing to do with the size of the original state space for
$N$ particles: it depends mostly on the smoothness of functions $g$
and $\phi$ rather than on $N$.
The second one is the
suppression of the expectation which might reduce the computational time by
a polynomial factor\footnote{The size of $\Pr_N(\SS)$ is the binomial 
  coefficient $(N{+}1{+}S,S)\sim_{N\to\infty} \frac{N^S}{S!}$} by replacing
the $|\Pr_N(\SS)|$ possible values of $M\toN_{t+1}$ by $1$.  The last
one  is
that the suppression of this expectation allows one  to carry the computation 
going forward rather than backward. This latter point  is particularly useful when the action set and the
time horizon are small.

\section{Conclusion and future work}

In this paper, we have shown how the mean field framework can be used
in an optimization context: the results known  for Markov chains can be
transposed almost unchanged to  Markov decision processes.
We further show that the convergence to the mean field limit in both
cases (Markovian and  Markovian with controlled variables) satisfies a
central limit theorem, providing  insight on the speed of
convergence.

We are currently investigating several extensions of these results.  First,
if one allows the actions to depend on the particles, it seems natural that
the limit behavior of such systems is the same as the limit behavior of
systems where the actions are random variables and that they both converge
to mean field system whose cost is averaged.  Another possible direction is
to consider stochastic systems where the event rate depends on $N$. In such
cases the deterministic limits are given by differential equations and the
speed of convergence can also be studied.

\bibliographystyle{plain}
\bibliography{bibFile}

\begin{thebibliography}{10}

\bibitem{EGEE}
{EGEE: Enabling Grids for E-sciencE}.

\bibitem{venkatBordenave}
V.~Anantharam and C.~Bordenave.
\newblock Optimal control of interacting particle systems.
\newblock {\em Private Communication}, 2008.

\bibitem{benaim:cmf}
M.~Bena{\i}m and J.Y. Le~Boudec.
\newblock {A Class Of Mean Field Interaction Models for Computer and
  Communication Systems}.
\newblock {\em To appear in Performance Evaluation}.

\bibitem{ArticleBerten.BG_PC07}
Vandy Berten and Bruno Gaujal.
\newblock Brokering strategies in computational grids using stochastic
  prediction models.
\newblock {\em Parallel Computing}, 2007.
\newblock Special Issue on Large Scale Grids.

\bibitem{InProceedingsBerten.BG_gbfbau_07}
Vandy Berten and Bruno Gaujal.
\newblock Grid brokering for batch allocation using indexes.
\newblock In {\em Euro-FGI NET-COOP}, Avignon, France, 2007. LNCS.

\bibitem{bordenave2007psi}
C.~Bordenave, D.~McDonald, and A.~Proutiere.
\newblock {A particle system in interaction with a rapidly varying environment:
  Mean field limits and applications}.
\newblock {\em Arxiv preprint math.PR/0701363}, 2007.

\bibitem{borkarStochasticApprox}
V.~Borkar.
\newblock {\em Stochastic Approximation: A Dynamical Systems Viewpoint}.
\newblock Cambridge University Press, 2008.

\bibitem{boudec2007gmf}
J.Y.L. Boudec, D.~McDonald, and J.~Mundinger.
\newblock {A Generic Mean Field Convergence Result for Systems of Interacting
  Objects}.
\newblock {\em QEST 2007.}, pages 3--18, 2007.

\bibitem{durrett1991pta}
R.~Durrett.
\newblock {\em {Probability: theory and examples}}.
\newblock Wadsworth \& Brooks/Cole, 1991.

\bibitem{Graham}
Carl Graham.
\newblock Chaoticity on path space for a queueing network with selection of the
  shortest queue among several.
\newblock {\em Journal of Applied Probability}, 37:198--211, 2000.

\bibitem{hochbaum1996aan}
D.S. Hochbaum.
\newblock {\em {Approximation algorithms for NP-hard problems}}.
\newblock PWS Publishing Co. Boston, MA, USA, 1996.

\bibitem{Mitrani}
Jennie Palmer and Isi Mitrani.
\newblock Optimal and heuristic policies for dynamic server allocation.
\newblock {\em Journal of Parallel and Distributed Computing},
  65(10):1204--1211, 2005.
\newblock Special issue: Design and Performance of Networks for Super-,
  Cluster-, and Grid-Computing (Part I).

\bibitem{Tsitsiklis}
Christos~H. Papadimitriou and John~N. Tsitsiklis.
\newblock The complexity of optimal queueing network control.
\newblock {\em Math. Oper. Res.}, 24:293--305, 1999.

\bibitem{puterman1994mdp}
M.L. Puterman.
\newblock {\em {Markov Decision Processes: Discrete Stochastic Dynamic
  Programming}}.
\newblock John Wiley \& Sons, Inc. New York, NY, USA, 1994.

\bibitem{Rolski}
T.~Rolski.
\newblock Comparison theorems for queues with dependent interarrival times.
\newblock In {\em Lecture Notes in Control and Information Sciences},
  volume~60, pages 42--71. Springer-Verlag, 1983.

\bibitem{Weber}
Richard~R. Weber and Gideon. Weiss.
\newblock On an index policy for restless bandits.
\newblock {\em Journal of Applied Probability}, 27:637--648, 1990.

\bibitem{Whittle}
P.~Whittle.
\newblock {\em A celebration of applied probability}, volume 25A, chapter
  Restless bandits: activity allocation in a changing world, pages 287--298.
\newblock J. Appl. Probab. Spec., j. gani edition, 1988.

\end{thebibliography}

\nocite{venkatBordenave}



\end{document}